\documentclass[12pt]{article}
\title{Birational geometry for nilpotent orbits}
\author{Yoshinori Namikawa}
\date{ }
\usepackage{amsmath,amssymb,amsthm, amscd}
\chardef\bslash=`\\

\def\0{{\mathcal O}}
\def\g{{\mathfrak g}}
\def\p{{\mathfrak p}}
\def\h{{\mathfrak h}}

\def\n{{\mathfrak n}}
\def\q{{\mathfrak q}}
\def\r{{\mathfrak r}}
\def\k{{\mathfrak k}}
\def\l{{\mathfrak l}}
\def\b{{\mathfrak b}}
\def\z{{\mathfrak z}}

\begin{document}
\maketitle

{\bf  \S 1.  Introduction}

The aim of this paper is to give an account of 
the birational point of view on nilpotent orbits 
in a complex simple Lie algebra. Let $\g$ be a 
complex simple Lie algebra and $G$ the adjoint 
group. An adjoint orbit $O$ in $\g$ is called a 
nilpotent orbit if $O$ consists of nilpotent 
elements of $\g$. The closure $\bar{O}$ of $O$ 
is then an affine variety with singularities. 
In general, $\bar{O}$ is not necessarily normal (see 
for example \cite{K-P} in this direction).  
In this paper we shall take its normalization 
$\tilde{O}$ and consider the birational 
geometry on its (partial) resolutions. 
Each variety $\tilde{O}$ has symplectic singularities. 
More precisely, the smooth locus $\tilde{O}_{reg}$ 
admits the Kostant-Kirillov 2-form $\omega$, which is $d$-closed 
and non-degenerate. Moreover, if we take a resolution 
$\mu: Y \to \tilde{O}$, then $\omega$ extends to a regular 2-form 
on $Y$. A resolution $\mu: Y \to \tilde{O}$ is called a 
crepant resolution if $K_Y = \mu^*K_{\tilde O}$.  
The nilpotent cone $N$ is defined to be the subset of $\g$ 
which consists of all nilpotent elements of $\g$. 
By definition $N$ is a disjoint union of all nilpotent orbits of 
$\g$. There is a largest nilpotent orbit $O_r$ and $N$ coincides 
with its closure. Moreover, $N$ is a normal variety. 
Let $B$ be a Borel subgroup of $G$ and let $T^*(G/B)$ be the 
cotangent bundle of the flag variety $G/B$. 
By using the Killing form of $\g$, one can identify 
$T^*(G/B)$ with a vector bundle $G \times^B [\b, \b]$ over $G/B$. 
Then there is a natural map 
$$ \nu: G \times^B [\b, \b] \to \g$$ defined by $[g,x] \to Ad_g(x)$. 
The image of $\nu$ coincides with $N$ and $\nu$ 
gives a resolution of $N$ (\cite{Sp}). We call $\nu$ the {\em Springer resolution} of $N$. 
Since $T^*(G/B)$ admits a canonical symplectic 2-form and it 
coincides with the pull-back of the Kostant-Kirillov 2-form 
on $O_r$, the Springer resolution is a crepant resolution. 
One can generalize this construction to a parabolic subgroup 
$Q$ of $G$. Let us start with the cotangent bundle $T^*(G/Q)$. 
Note that $T^*(G/Q)$ is identified with $G \times^Q \n(\q)$ where 
$\n(\q)$ is the nil-radical of $\q$. 
In a similar way to the above, we have a map 
$$\nu: T^*(G/Q) \to \g,$$ whose image is the closure of 
a nilpotent orbit $O$. In general, $\nu$ is not birational onto 
its image, but a generically finite projective morphism (see Example (2.5.4) 
for a non-birational Springer map). 
When $\nu$ gives a resolution of $\bar{O}$, we call $\nu$ the 
Springer resolution of $\bar{O}$. In this case, the Stein factorization 
$$T^*(G/Q) \stackrel{\nu^n}\to \tilde{O} \to \bar{O}$$ gives 
a crepant resolution of $\tilde{O}$.   
B. Fu [Fu 1] proved the following. 

{\bf Theorem}. {\em Let $O$ be a nilpotent orbit of $\g$ and assume that 
$\tilde{O}$ admits a crepant resolution. Then it coincides with a Springer 
resolution. More exactly, there is a parabolic subgroup $Q$ of $G$ such that 
$\nu^n$ is the given crepant resolution.} 
\vspace{0.2cm}

However there still remained interesting problems. At first, there actually exists 
a nilpotent orbit which has no crepant resolutions.
Secondly, even if $\tilde{O}$ has a crepant resolution, it is not unique, that is, 
the choice of $Q$ is not unique even up to conjugacy class.  Our purpose is to survey  
complete answers (cf. \cite{Na 1}, \cite{Na 2} , \cite{Na 3} and \cite{Fu 2}) to these problems.  

A substitute for a crepant resolution is a {\em {\bf Q}-factorial terminalization}. 
A birational projective morphism $\mu: Y \to \tilde{O}$ is a {\bf Q}-factorial 
terminalization if $Y$ has only {\bf Q}-factorial terminal singularities and 
$K_Y = \mu^*K_{\tilde O}$. The exsistence of a {\bf Q}-factorial terminalization 
is established by Birkar, Cascini, Hacon and McKernan [BCHM]. 
But, we shall give here more concrete forms of {\bf Q}-factorial terminalization. 
A hint is already in the work of Lustzig and Spaltenstein [LT].    
They introduced the notion of an {\em induced orbit}. 
Let us start with a parabolic subgroup $Q$ of $G$ and its Levi factor $L(Q)$. 
Let $O' \subset \l(\q)$ be a nilpotent orbit with respect to the adjoint 
$L(Q)$-action. Then one can make an associated bundle $G \times^Q (\n(\q) + \bar{O'})$ 
and define a map $$\nu: G \times^Q (\n(\q) + \bar{O'}) \to \g$$ by 
$\nu ([g,x]) = Ad_g(x)$. Since this is a $G$-equivariant closed map, its image is the 
closure of a nilpotent orbit $O$ of $\g$. Then we say that $O$ is induced from 
$O'$ and write $O = \mathrm{Ind}^{\g}_{\l(\q)}(O')$. The map 
$\nu$ is called the generalized Springer map. The generalized Springer map $\nu$ is a 
generically finite projective morphism. But if $\nu$ is birational onto its image, 
then the Stein factorization 
$$ G \times^Q (\n(\q) + \bar{O'}) \stackrel{\nu^n}\to \tilde{O} \to \bar{O}$$ 
gives a partial resolution of $\tilde{O}$.     
Now one can prove: 

{\bf Theorem (2.6.2)}. {\em Let $O$ be a nilpotent orbit of a complex simple Lie 
algebra $\g$. Then there are a parabolic subalgebra $\q$ of $\g$ and a nilpotent orbit $O'$ of 
$\l(\q)$ such that the following holds}: 
\vspace{0.15cm}
 
(1) $O = \mathrm{Ind}^{\g}_{\l(\q)}(O')$. 

(2) {\em $\nu^n$ gives a {\bf Q}-factorial terminalization of $\tilde{O}$.} 
\vspace{0.15cm} 

In order to look for other {\bf Q}-factorial terminalizations of $\tilde{O}$, 
we introduce a flat deformation of $G \times^Q (\n(\q) + \tilde{O'})$. 
For simplicity we put $\l := \l(\q)$ and let $L$ be the corresponding Levi subgroup. 
Let $\r(\q)$ be the solvable radical of $\q$ and consider the variety 
$G \times^Q(\r(\q) + \bar{O'})$. Its normalization $X_{\q,O'}$ is isomorphic 
to $G \times^Q(\r(\q) + \tilde{O'})$. Let $\k$ be the center of $\l$. 
In (3.3) we shall define a map $$X_{\q,O'} \to \k$$ whose central fiber 
$X_{\q,O',0}$ is $G \times^Q (\n(\q) + \tilde{O'})$. This map factorizes as 
$$X_{\q,O'} \stackrel{\mu_{\q}}\to \mathrm{Spec}\; \Gamma (X_{\q,O'}, \mathcal{O}_{X_{\q,O'}}) 
\to \k.$$  Put 
$$Y_{\l,O'} := \mathrm{Spec}\; \Gamma (X_{\q,O'}, \mathcal{O}_{X_{\q,O'}}).$$ 
An important fact is that $Y_{\l,O'}$ depends only on $\l$ and $O'$. 
Moreover its central fiber $Y_{\l,O',0}$ is isomorphic to $\tilde{O}$. 
Define 
$$ \mathcal{S}(\l) := 
\{\mathrm{parabolic}\;  \mathrm{subalgebras}\;  \q' \; \mathrm{of} \; \g ;\; \l(\q') = \l \}.$$  
We can define $X_{\q',O'}$ for each $q' \in \mathcal{S}(\l)$. 
We also have a map $$\mu_{\q'}: X_{\q',O'} \to Y_{\l,O'}.$$
The map $\mu_{\q'}$ is a crepant birational morphism. Moreover, $\mu_{\q',t}$ 
is an isomorphism for $t \in \k^{reg}$; hence $\mu_{\q'}$ is an isomorphism in 
codimension one. Define $$M(L) := \mathrm{Hom}_{alg.gp}(L, \mathbf{C}^*)$$ 
and put $M(L)_{\mathbf R} := M(L) \otimes \mathbf{R}$. Then 2-nd cohomology 
groups $H^2(X_{\q',O'}, \mathbf{R})$ are naturally identified with $M(L)_{\mathbf R}$. 
By these identifications the nef cones $\overline{\mathrm{Amp}}(\mu_{\q',O'})$ are 
regarded as the cones in $M(L)_{\mathbf R}$. 
\vspace{0.2cm}

{\bf Theorem (3.5.1)} {\em For $\q' \in \mathcal{S}(\l)$, the birational map 
$\mu_{\q'}: X_{\q',O'} \to Y_{\l,O'}$ is a {\bf Q}-factorial terminalization and 
is an isomorphism in codimension one.  
Any {\bf Q}-factorial terminalization of $Y_{\l, O'}$ is obtained in this way. 
If $\q_1 \ne \q_2$, then $\mu_{\q_1}$ and $\mu_{\q_2}$ give different 
{\bf Q}-factorial terminalizations. Moreover,} 
$$M(L)_{\mathbf R} = \cup_{\q' \in \mathcal{S}(\l)} \overline{\mathrm{Amp}}(\mu_{\q'}).$$
  
Two elements of $\mathcal{S}(\l)$ are connected by a sequence of the operations called {\em 
twists} (cf. (3.1)). Corresponding to a twist $\q_1 \leadsto \q_2$, we have a flop 
$$X_{\q_1,O'} \to Z \leftarrow X_{\q_2,O'}.$$ So any two {\bf Q}-factorial terminalizations of $\tilde{O}$ 
are connected by a suequence of certain flops. Now let us look at the central fibers 
$X_{\q',O',0}$ of $X_{\q',O'} \to \k$.  The diagram $$X_{\q_1,O',0} \to Z_0 \leftarrow X_{\q_2,O',0}$$ is not 
necessarily  a flop. 
Twists are divided into those of the first kind and those of the second kind. 
If the twist $\q_1 \leadsto \q_2$ is of the first kind, then it induces a flop between  
$X_{\q_1,O',0}$ and $X_{\q_2,O',0}$. These flops are completely classified and 
we call them Mukai flops (3.2.1). If it is of the second kind, the maps 
$X_{\q_i,O',0} \to Z_0$ ($i = 1,2$) are both divisorial birational maps.     
Define $\mathcal{S}^1(\l)$ to be the subset of $\mathcal{S}(\l)$ consisting of the 
parabolic subalgebras $\q'$ obtained from $\q$ by a finite succession of the twists of 
the first kind. 
Note that the restriction map $H^2(X_{\q',O'}, \mathbf{R}) 
\to H^2(X_{\q',O',0}, \mathbf{R})$ is an isomorphism and $\overline{\mathrm{Amp}}(\mu_{\q'})$ 
is mapped onto $\overline{\mathrm{Amp}}(\mu_{\q',0})$. 
\vspace{0.2cm}

{\bf Theorem (3.5.4)}. {\em There is a one-to-one correspondence between the set of 
{\bf Q}-factorial terminalizations of $\tilde{O}$ and $\mathcal{S}^1(\l)$. 
In other words, every {\bf Q}-factorial terminalization of $\tilde{O}$ is obtained as  
$\mu_{\q',0}: X_{\q,O',0} \to \bar{O}$ for $\q' \in 
\mathcal{S}^1(\l)$. Two different {\bf Q}-factorial terminalizations of 
$\tilde{O}$ are connected by a sequence of Mukai flops. 
Moreover $$\overline{\mathrm{Mov}}(\mu_{\q,0}) = \cup_{\q' \in 
\mathcal{S}^1(\l)} \overline{\mathrm{Amp}}(\mu_{\q',0}),$$ where 
$\overline{\mathrm{Mov}}(\mu_{\q,0})$ is the movable cone for $\mu_{\q,0}$ 
(cf. (3.4.2)). } 
\vspace{0.2cm}

A direct approach to Theorem (3.5.4) usually needs the classification 
of the generalized Springer maps which are isomorphisms in codimension 
one. But our approach using Theorem (3.5.1) does not need this and 
Mukai flops appear in a very natural way. 

Let $W$ be the Weyl group of $\g$ and let $N_W(L)$ be the subgroup 
of $W$ which normalizes $L$. Then the quotient group 
$$ W' := N_W(L) /W(L) $$ naturally acts on $M(L)_{\mathbf R}$. 
The interior ${\mathrm{Mov}}(\mu_{\q,0})$ of the movable cone can be characterized 
as a fundamental domain for this action (Theorem (3.6.1)). 
The group $W'$ was extensively studied in [Ho].  
As explained above, the deformation $X_{\q,O'} \to \k$ of 
$G \times^Q (\n(\q) + \tilde{O}')$ played an important role 
to study the birational geometry for $\tilde{O}$. But this is 
not merely a flat deformation of $G \times^Q (\n(\q) + \tilde{O}')$.   
In fact, $G \times^Q (\n(\q) + \tilde{O}')$ admits a symplectic 
2-form on its regular locus. This symplectic 2-form induces a 
Poisson structure of the regular part; moreover, it uniquely 
extends to a Poisson structure of $G \times^Q (\n(\q) + \tilde{O}')$. 
One can introduce the notion of a Poisson deformation (cf. \S 4), 
and $X_{\q,O'} \to \k$ turns out to be a Poisson deformation of 
$G \times^Q (\n(\q) + \tilde{O}')$. On the other hand, since $\tilde{O}$ 
has symplectic singularities, $\tilde{O}$ also admits a natural Poisson structure. 
One can construct a flat deformation of $\tilde{O}$ as follows. 
Let $G \cdot (\r(\q) + \bar{O}') \subset \g$ denote the $G$-orbit of $\r(\q) + \bar{O}'$ 
by the adjoint $G$-action. By using the adjoint quotient map $\g \to \h/W$, 
we get a map $\chi: G \cdot (\r(\q) + \bar{O}') \to \h/W$. The image 
of  $\chi$ is not necessarily normal, but its normalization coincides with 
$\k(\q)/W'$. Let $G \cdot (\r(\q) + \bar{O}')^n$ be the normalization 
of  $G \cdot (\r(\q) + \bar{O}')$. Then $\chi$ induces a map 
$$ G \cdot (\r(\q) + \bar{O}')^n \to \k(\q)/W'.$$  
One can check that this is a flat map and its central fiber is $\tilde{O}$. 
Moreover, this is a Poisson deformation of $\tilde{O}$. 
The two Poisson deformations are combined together by the 
Brieskorn-Slodowy diagram 

\begin{equation} 
\begin{CD} 
X_{\q,O'} @>>> G \cdot(\r(\q) + \bar{O}')^n \\ 
@VVV @VVV \\ 
\k(\q) @>>> \k(\q)/W'     
\end{CD} 
\end{equation}

Our last theorem (= Theorem (4.5)) claims that it gives formally 
universal Poisson deformations of 
$G \times^Q (\n(\q) + \tilde{O}')$ and $\tilde{O}$.  
 
Finally we shall explain the contents of this paper. 
The first part of \S 2 is an introduction to nilpotent orbits 
and related resolutions. Many concrete examples are desrcibed 
in terms of flags; I believe that they would motivate the following 
abstract arguments. In the final part of \S 2, we give a rough sketch 
of the proof of Theorem (2.6.2) in the classical cases. 
The readers can find a proof in [Fu 2] when $\g$ is exceptional.    
The idea of most arguments in \S 3 comes from [Na 2]. But all statements 
are generalized so that one can apply them to generalized Springer maps. 
\S 4 is concerned with a Poisson deformation. We quickly review the 
notions of Poisson structures and Poisson deformations. After that, we 
will give a rough sketch of Theorem (4.5) mentioned above. The results 
of \S 4 have been already treated in [Na 5] when $\tilde{O}$ has a 
crepant resolution.  
\vspace{0.2cm} 

{\bf Notations}. Let $G$ be an algebraic group over ${\mathbf C}$ and 
$P$ a closed subgroup of $G$. If $V$ is a variety with a $P$-action, 
then we denote by $G \times^P V$ the associated fiber bundle over $G/P$ 
with a typical fiber $V$. More exactly, $G \times^P V$ is defined as 
the quotient of $G \times P$ by an equivalence relation $\sim$, where 
$(g,x) \sim (g',x')$ if there is an element $p \in P$ such that 
$g' = gp$ and $x' = p^{-1}\cdot x$.    
\vspace{0.2cm}

{\bf  \S 2. Nilpotent orbits and symplectic singularities}

(2.1) Let $G$ be a semi-simple algebraic group over the complex number field 
$\mathbf{C}$ and let $\g$ be its Lie algebra. An orbit $O$ of the adjoint 
action $Ad: G \to \mathrm{Aut}(\g)$ is called an adjoint orbit. Moreover, 
if $O$ consists of nilpotent elements (resp. semi-simple elements), then $O$ is called a 
{\em nilpotent orbit} (resp. {\em semi-simple orbit}). 
The tangent space $T_{\alpha}O$ of an adjoint orbit $O$ 
at $\alpha$ is identified with 
$$[\alpha, \g] := \{[\alpha, x]; x \in \g\}.$$ 
Since $\g$ is semi-simple, the Kostant-Kirillov form 
$$ k: \g \times \g \to \mathbf{C}$$ is a non-degenerate symmetric 
form. We define a skew-symmetric form 
$$ \omega_{\alpha}: T_{\alpha}O \times T_{\alpha}O \to \mathbf{C}$$ by 
$$\omega_{\alpha}([\alpha, x], [\alpha, y]) := k(\alpha, [x,y]).$$ 
  This is well-defined and non-degenerate because if $[\alpha, x] = 0$, 
then $k(\alpha, [x,y]) = k([\alpha, x], y) = 0$. If $\alpha$ runs through 
all elements of $O$, the 2-form $\omega := \{\omega_{\alpha}\}$ is a 
{\em d-closed} form on $O$. In particular, $O$ is a smooth algebraic variety of 
even dimension. The symplectic form $\omega$ is called the 
Kostant-Kirillov 2-form. A semi-simple orbit is a closed subvariety of $\g$. 
But, a nilpotent orbit $O$ is not closed in $\g$ except when $O = \{0\}$. 
If we take the closure  $\bar{O}$ of $O$, it is an affine variety with singularities. 
Note that $\bar{O}$ is not necessarily normal. We denote by $\tilde{O}$ its 
normalization.    
\vspace{0.2cm}

(2.2) {\em Nilpotent orbits in a classical Lie algebra}: 
Let $sl(n)$ be the Lie algebra consisting of $n \times n$ matrices $A$ with 
$tr(A) = 0$. Define $$so(n) := \{A \in sl(n); A^tJ + JA = 0\}, $$ where    
$$ J =  
\begin{pmatrix} 
 & & & & & 1\\
 & & & & . &\\
 & & & . & &\\
 & & . & & &\\
 & . & & & &\\  
1 & & & & & 
\end{pmatrix}, 
$$  and $A^t$ is the transposed matrix of $A$. 
Similarly, define $$sp(2n) := \{A \in sl(2n); A^tJ' + J'A = 0\}, $$ where 
$$ J' =  
\begin{pmatrix} 
& & & & & & & 1\\
& & & & & & . &\\
& & & & & . & &\\

& & & & 1 & & &\\ 
& & & -1& & & & \\
& & . & & & & &\\
& . & & & & & &\\
-1 & & & & & & &
\end{pmatrix}.
$$  

If $\g$ is of type $A_{n-1}$, then $\g = sl(n)$. If $\g$ is 
of type $B_n$, then $\g = so(2n+1)$.  If $\g$ is of type 
$C_n$, then $\g = sp(2n)$. Finally, if $\g$ 
is of type $D_n$, then $\g = so(2n)$.     
One can associate a Jordan type $\mathbf{d}$ to each nilpotent orbit 
of $\g$. Here $\mathbf{d} := [d_1, d_2, ..., d_k]$ is a partition of 
$n := \dim V$. Namely, $d_i$ are positive integers such that $d_1 \geq d_2 \geq 
... \geq d_k$ and $\Sigma d_i = n$.        
Another way of writing $\mathbf{d}$ is 
$[d_1^{s_1}, ..., d_k^{s_k}]$ with $d_1 > d_2 ... > d_k >0$. 
Here $d_i^{s_i}$ is an $s_i$ times $d_i$'s: 
$d_i, d_i, ..., d_i$.   
The partition  ${\mathbf d}$ corresponds to a Young diagram. 
For example, $[5,  4^2,  1]$ corresponds to 

\begin{picture}(100, 100)(0, 0)
\put(00,  80){\line(1, 0){100}}
\put(00,  60){\line(1, 0){100}}
\put(00,  40){\line(1, 0){80}}
\put(00,  20){\line(1, 0){80}}
\put(00,  00){\line(1, 0){20}}
\put(00,  00){\line(0, 1){80}}
\put(20,  00){\line(0, 1){80}}
\put(40,  20){\line(0, 1){60}}
\put(60,  20){\line(0, 1){60}}
\put(80,  20){\line(0, 1){60}}
\put(100,  60){\line(0, 1){20}}
\end{picture}

When an   
integer $e$ appears in the partition $\bf{d}$, we say 
that $e$ is a {\em member} of $\bf{d}$. We call $\bf{d}$ 
{\em very even} when $\bf{d}$ consists with only even 
members, each having even multiplicity. 
  
The following result can be found, for example, in 
[C-M, \S 5].     

{\bf Proposition (2.2.1)}
{\em   Let $\mathcal{N}o(\g)$ be the set of nilpotent 
orbits of $\g$.} 
\vspace{0.12cm}

(1)($A_{n-1}$): {\em When $\g = \mathfrak{sl}(n)$, there is a 
bijection between $\mathcal{N}o(\g)$ and 
the set of partitions $\bf{d}$ of $n$. }  
\vspace{0.12cm} 
 
(2)($B_n$): {\em When $\g = \mathfrak{so}(2n+1)$, there is a 
bijection between $\mathcal{N}o(\g)$ and 
the set of the partitions $\bf{d}$ of $2n+1$ 
for which all even members have even multiplicities.}    
\vspace{0.12cm}

(3)($C_n$): {\em When $\g = \mathfrak{sp}(2n)$, there is a 
bijection between $\mathcal{N}o(\g)$ and the 
set of the partitions $\bf{d}$ of $2n$ for which all odd members 
have even multiplicities.}
\vspace{0.12cm}

(4)($D_n$): {\em When $\g = \mathfrak{so}(2n)$, there is a 
surjection $f$ from $\mathcal{N}o(\g)$ to the set 
of the partitions $\bf{d}$ of $2n$ for which  
all even members have even multiplicities. For a partition $\bf{d}$ which is 
not very even, $f^{-1}(\bf{d})$ consists of exactly one orbit, 
but, for very even $\bf{d}$, $f^{-1}(\bf{d})$ consists of exactly 
two different orbits. }  \vspace{0.2cm}

We introduce a partial order in the set of the  
partitions of (the same number):  
for two partitions $\mathbf{d}$ and $\mathbf{f}$,  
$\mathbf{d} \geq \mathbf{f}$ if $\Sigma_{i \leq k}d_i 
\geq \Sigma_{i \leq k}f_i$ for all $k \geq 1$.  
On the other hand, for two nilpotent orbits $\mathcal{O}$ 
and $\mathcal{O}'$ in  $\g$, we write $\mathcal{O} \geq {\mathcal O}'$ 
if $\mathcal{O}' \subset \bar{\mathcal{O}}$. 
Then, $\mathcal{O}_{\mathbf d} \geq \mathcal{O}_{\mathbf f}$ if and 
only if $\mathbf{d} \geq \mathbf{f}$. 

{\bf Remark (2.2.2)}. In order to classify nilpotent orbits in  
simple Lie algebras including those of exceptional type, we need 
a different method. Dynkin [D] associates a weighted Dynkin diagram to 
each nilpotent orbit. This correspondence is an injection, but is not 
surjective. Bala and Carter [B-C] determined which weighted Dynkin 
diagrams come from nilpotent orbits and completed the classification of 
nilpotent orbits in all simple Lie algebras. For details, see [B-C] and [C-M].   
\vspace{0.2cm}

(2.3) {\em Jacobson-Morozov resolution of $\bar{O}$}: 
Let $O$ be a nilpotent orbit of a complex semi-simple Lie algebra $\g$. 
Fix an element $x \in O$.   
By the Jacobson-Morozov theorem (cf. [C-M, 3.2]) one can find 
a semi-simple element $h \in \g$, and a 
nilpotent element $y \in \g$ in such a way that 
$[h,x] = 2x$, $[h,y] = -2y$ and $[x,y] = h$. 
For $i \in \mathbf{Z}$, let 
$$\g_i := \{z \in \g\; [h,z] = iz \}.$$ 
Then one can write 
$$ \g = \oplus_{i \in \mathbf{Z}}\g_i. $$ 
Let $\h$ be a Cartan subalgebra of $\g$ with 
$h \in \h$.  Let $\Phi$ be the corresponding 
root system and let $\Delta$ be a base of simple roots 
such that $h$ is $\Delta$-dominant, i.e. $\alpha(h) 
\geq 0$ for all $\alpha \in \Delta$. 
In this situation, 
$$\alpha(h) \in \{0,1,2\}.$$ 
The weighted Dynkin diagram of $\mathcal{O}_x$ is 
the Dynkin diagram of $\g$ where each vertex $\alpha$ is 
labeled with $\alpha(h)$. A nilpotent orbit $\mathcal{O}_x$ 
is completely determined by its weighted Dynkin diagram. 
A Jacobson-Morozov parabolic subalgebra for 
$x$ is the parabolic subalgebra $\p$ defined by 
$$\p:= \oplus_{i \geq 0}\g_i.$$ 
Let $P$ be the parabolic subgroup of $G$ determined by  $\p$. 
We put $$ \n_2 := \oplus_{i \geq 2}\g_i. $$ 
Then $\n_2$ is an ideal of $\p$; hence,  
$P$ has the adjoint action on $\n_2$. Let us consider the   
vector bundle  $G \times^P \n_2$ over $G/P$ and the map 
$$\mu: G \times^P \n_2 \to \g$$ defined by $\mu([g,z]) := 
Ad_g(z)$. Then the image of $\mu$ coincides with the 
closure $\bar{O}$ of $O$ and 
$\mu$ gives a resolution of $\bar{O}$. 
We call $\mu$ the {\em Jacobson-Morozov resolution} 
of $\bar{O}$. The construction of $\mu$ depends on the 
choices of $x \in O$ and the $sl(2)$-triple $\{x,y,h\}$. 
But, for any nilpotent elements $x$ and $x'$ of $O$, 
two $sl(2)$-triplet $\{x,y,h\}$ and $\{x', y', h'\}$ are 
conjugate to each other by an element of $G$ (cf. [C-M, 3.2]). In this sense,  
the Jacobson-Morozov resolution is unique. But, 
the Jacobson-Morozov resolution is {\em not} a crepant 
resolution in general.   
\vspace{0.2cm}

{\bf Definition (2.4)}.  Let $X$ be a normal variety defined over $\mathbf{C}$. 
Assume that the regular locus $X_{reg}$ admits a symplectic 2-form $\omega$. 
Then $(X, \omega)$ (or $X$) has {\em symplectic singularities} if there is a  
resolution $\mu: Y \to X$ such that the 2-form $\omega$ on $\mu^{-1}(X_{reg})$  
extends to a regular 2-form on $Y$. 
\vspace{0.2cm}

Remark that if a particular resolution $\mu: Y \to X$ has the extension property 
explained above, then all resolutions of $X$ actually have the property.   
The following proposition is due to Hinich and Panyushev [Hi],[Pa]:  

{\bf Proposition (2.4.1)}. {\em For the Jacobson-Morozov resolution $\mu: 
G \times^{P}\n_2 \to \bar{O}$, the Kostant-Kirillov 2-form on $O$ extends 
to a regular 2-form on $G \times^{P}\n_2$. In particular, the normalization 
$\tilde{O}$ of $\bar{O}$ has symplectic singularities.}   
\vspace{0.2cm}

(2.5) ({\em Induced orbits}):  

(2.5.1) Let $G$ and $\g$ be the same as in (2.1). 
Let $Q$ be a parabolic subgroup of $G$ and let $\q$ be its Lie 
algebra with Levi decomposition $\q = \l \oplus \n$. Here $n$ is the 
nil-radical of $\q$ and $\l$ is a Levi-part of $\q$. 
Fix a nilpotent orbit $O'$ in $\l$. Then there is a unique 
nilpotent orbit $O$ in $\g$ meeting $n + O'$ in 
an open dense subset ([L-S]). Such an orbit $O$ is called 
the nilpotent orbit induced from $O'$ and we write $$O = 
\mathrm{Ind}^{\g}_{\l}(O').$$  
Note that when $O' = 0$, $O$ is called the  
{\em Richardson orbit} for $Q$. 
Since the adjoint action of 
$Q$ on $\q$ stabilizes $n + \bar{O'}$, one can consider 
the variety $G \times^Q (n + \bar{O'})$.   
There is a map 
$$ \nu: 
G \times^Q (\n + \bar{O'}) \to \bar{O}$$ defined 
by $\nu ([g,z]) := Ad_g(z)$. 
Since $\mathrm{Codim}_{\l}(O') = 
\mathrm{Codim}_{\g}(O)$ (cf. [C-M], Prop. 7.1.4), 
$\nu$ is a generically finite dominating map. Moreover, $\nu$ is factorized as 
$$G \times^Q (\n + \bar{O'}) \to G/Q \times \bar{O} 
\to \bar{O}$$ where the first map is a closed embedding 
and the second map is the 2-nd projection; this implies that 
$\nu$ is a projective map. In the remainder, we call $\nu$ the 
{\em generalized Springer map} for ($Q$, $O'$). When 
$O' = \{0\}$, we often call $\mu$ the Springer map. 
Let $\tilde{O'}$ be the normalization of $\bar{O'}$. 
Then the normalization of $G \times^Q (\n + \bar{O'})$ coincides 
with $G \times^Q (\n + \tilde{O'})$. The generalized Springer map 
$\nu$ induces a map $$\nu^n: G \times^Q (\n + \tilde{O'}) \to 
\tilde{O}.$$ We call $\nu^n$ the {\em normalized map} of  $\nu$. 

(2.5.2) Assume that there are a parabolic subgroup $Q_L$ of $L$ and 
a nilpotent orbit $\mathcal{O}''$ in the Levi subalgebra $\l(Q_L)$ such that 
$O'$ is the nilpotent orbit induced from $(Q_L, O'')$. 
Then there is a parabolic subgroup $Q'$ of $G$ such that $Q' \subset Q$, 
$\l(Q') = \l(Q_L)$ and $\mathcal{O}$ is the nilpotent orbit induced from 
$(Q', O'')$. 
If we put $\l' := \l(Q')$, then this can be written as 
$$\mathrm{Ind}^{\g}_{\l}(\mathrm{Ind}^{\l}_{\l'}(O'')) = 
\mathrm{Ind}^{\g}_{\l'}(O'').$$
The generalized Springer map $\nu'$ for $(Q', O'')$ 
is factorized as 
$$ G \times^{Q'}(\n' + \bar{O''}) \to G \times^Q(\n + \bar{O'}) \to \bar{O}.$$ 

{\bf Example (2.5.3)}.  Let $G := SL(n, \mathbf{C})$ and $\g := sl(n)$. 
Then $G$ acts naturally on $\mathbf{C}^n$ and any parabolic subgroup 
$Q$ of $G$ is given as the subgroup of stabilizers of a flag (= a sequence 
of vector subspaces):  
$$V_1 \subset V_2 \subset ... \subset V_l = \mathbf{C}^n.$$ 
Put $q_i := \dim V_i - \dim V_{i-1}$ and $(q_1, q_2, ..., q_l)$ is 
called the type of $Q$. If two parabolic subgroups of $G$ have the same 
type, then they are conjugate to each other. 
The set of all diagonal matrices in $\g$ forms a Cartan subalgebra $\h$.  
There is a unique Levi decomposition $\q = \l \oplus \n$ such that 
$\h \subset \l$.  In our case, 
$$ \l = \{ \left( \begin{array}{cccc}
A_{q_1} & 0 & \cdots & 0 \\ 
0 & A_{q_2} & 0 & \cdots \\ 
\cdots & \cdots & \cdots & \cdots \\
0 & \cdots & 0 & A_{q_l} 
\end{array}\right) \in \g \; \vert  \;  A_{q_i} : q_i \times q_i \; \mathrm{matrix} \}.$$ 
We take $O': = \{0\}$ as a nilpotent orbit in $\l$.  
Rearrange $q_i's$ in such a way that $q_{\sigma(1)} \geq q_{\sigma(2)} \geq ... \geq q_{\sigma(l)}$ 
by using a suitable permutation $\sigma \in S_l$. 
Then $\mathbf{q} := (q_{\sigma(1)}, q_{\sigma(2)}, ..., q_{\sigma(l)})$ is a partition of $n$. As in (2.4),  
we associate a Young diagram to $\mathbf{q}$. Let $d_i$ be the length of the 
$i$-th column of the Young diagram. The dual partition $\mathbf{q}^t$ of $\mathbf{q}$ is 
defined as $\mathbf{q}^t := [d_1, d_2, ..., d_s]$.   We shall prove that 
$\mathrm{Ind}^{\g}_{\l}(\mathcal{O}') \subset \g$ is the nilpotent orbit with Jordan 
type $\mathbf{q}^t$.  Take a basis $e_1$, $e_2$, ..., $e_n$ of $\mathbf{C}^n$ in 
such a way that $$V_i = \mathbf{C}<e_1, e_2, ..., e_{\Sigma_{1 \le j \le i}q_j}>$$ for 
all $i$. Define $$W_1 := \{e_1, ..., e_{q_1}\},$$ 
$$W_2 := \{e_{q_1 + 1}, ..., e_{q_1 + q_2}\},$$ 
$$ ...... $$ $$W_l := \{e_{(\Sigma_{1 \le k \le l-1}q_k) +1}, ..., e_{\Sigma_{1 \le k \le l}q_k}\}.$$  
Then $V_i $ is spanned by $W_1 \cup ... \cup W_i$.  Take the 1-st vectors from $W_i$'s and 
form a set $E_1$ consisting of them. Namely $$E_1 = \{e_1, e_{q_1 + 1}, ..., e_{q_1 + ... + q_{l-1} +1}\}.$$ 
Next take the 2-nd vectors from $W_i$'s.  If $q_i = 1$, 
then there is no 2-nd vector in $W_i$. In this case, we take no vectors 
from this $W_i$. Let $E_2$ be the set consisting of such vectors.    
Similarly we define $E_3$, ... Note that $E_i$ has exactly $d_i$ elements.         
Let us consider the nilpotent endomorphism $x_i$ of $\mathbf{C}<E_i>$ corresponding 
to the Jordan matrix $J_{d_i}$ of size $d_i$:   
$$ J_{d_i} :=  \left( \begin{array}{ccccc} 
0 & 1 & 0 & \cdots & \cdots \\ 
0 & 0 & 1 & 0 & \cdots \\  
\cdots & \cdots & \cdots & \cdots & \cdots \\ 
0 & 0 & \cdots & 0 & 1 \\ 
0 & 0 & \cdots & 0 & 0 
\end{array} \right).  $$ 
We define a nilpotent endomorphism $x$ of $\mathbf{C}^n$ 
by $x := \oplus x_i$.  By the construction, $x$ has Jordan type 
$[d_1, d_2, ..., d_s]$ and $x(V_i) \subset V_{i-1}$ for 
each $i$; hence $x \in \n$. Let us consider the Springer map 
$$\nu: G \times^Q \n \to \g.$$ Since $x \in \n$,  the image of $\mu$ 
contains $x$. Let $O_x$ be the nilpotent orbit containing $x$. In 
order to prove that $\mathrm{Im}(\mu) = \bar{O}_x$, it suffices 
to prove that $\dim O_x = \dim (G \times^Q \n)$. We put 
$$ \g^x := \{z \in \g; [x,z] = 0\}.$$ Note that $\dim O_x = \dim \g - 
\dim \g^x$. By using this fact, one can check that  
$$\dim O_x = n(n+1) - 2\Sigma_{1 \le i \le s}id_i.$$ On the other hand, 
one has $$\dim (G/Q) = n(n+1)/2 - \Sigma_{1 \le i \le s}id_i.$$ 
Since $\dim (G \times^Q \n) = 2\dim (G/Q)$, we have the desired 
result. Finally we shall check that $\nu$ is a birational map. Since 
$\nu$ is a $G$-equivariant map, we only have to show that 
$\nu^{-1}(y)$ consists of exactly one point for a particular nilpotent 
element $y$ with Jordan type $[d_1, d_2, ..., d_s]$. 
We put 
$$ y = \left( \begin{array}{cccc}
J_{d_1} & 0 & \cdots & 0 \\ 
0 & J_{d_2} & 0 & \cdots \\ 
\cdots & \cdots & \cdots & \cdots \\
0 & \cdots & 0 & J_{d_s} 
\end{array}\right). $$  
Assume that $\nu([g,z]) = y$. Then $z$ is uniquely determined by $[g]$ as 
$z := Ad_{g^{-1}}(y)$. Therefore, $$\nu^{-1}(y) =  
\{[g] \in G/Q; y \in Ad_g(\n)\}.$$  
Note that $G/Q$ is naturally identified with the set of parabolic subgroups of $G$ 
which are conjugate to $Q$ by $[g] \to gQg^{-1}$. Moreover, $Ad_g(\n) 
= \n(gQg^{-1})$. As a consequence, $\nu^{-1}(y)$ is identified with the 
set of parabolic subgroups $Q'$ of $G$ such that $Q'$ is conjugate to 
$Q$ and $y \in \n(Q')$. In terms of flags, $Q'$ corresponds to a 
flag $V_1 \subset V_2 \subset ... \subset V_l = \mathbf{C}^n$ 
of type $(q_1,q_2, ..., q_l)$ so that $y(V_i) \subset V_{i-1}$ for all $i$.      
Let us consider the Young diagram corresponding to $[d_1, ..., d_s]$. 
We fill up each box on the 1-st row by $e_1$, ..., $e_{d_1}$ from left to right. 
Next fill up each box on the 2-nd row by $e_{d_1+1}$, ..., $e_{d_2}$ from 
left to right, and so on.  For example, when $(q_1, q_2, q_3, q_4, q_5) = (3,4,3,3,1)$ and 
$[d_1, d_2, d_3, d_4] = [5, 4, 4, 1]$, we have the following tablaux:    

\begin{picture}(100, 100)(0, 0)
\put(00,  80){\line(1, 0){100}}
\put(00,  60){\line(1, 0){100}}
\put(00,  40){\line(1, 0){80}}
\put(00,  20){\line(1, 0){80}}
\put(00,  00){\line(1, 0){20}}
\put(00,  00){\line(0, 1){80}}
\put(20,  00){\line(0, 1){80}}
\put(40,  20){\line(0, 1){60}}
\put(60,  20){\line(0, 1){60}}
\put(80,  20){\line(0, 1){60}}
\put(100,  60){\line(0, 1){20}}
\put(05, 65){$e_1$}
\put(25, 65){$e_2$}
\put(45, 65){$e_3$}
\put(65, 65){$e_4$}
\put(85, 65){$e_5$}
\put(05, 45){$e_6$}
\put(25, 45){$e_7$}
\put(45, 45){$e_8$}
\put(65, 45){$e_9$}
\put(05, 25){$e_{10}$}
\put(25, 25){$e_{11}$}
\put(45, 25){$e_{12}$}
\put(65, 25){$e_{13}$}
\put(05, 05){$e_{14}$}
\end{picture}

Let us consider the $q_1$ boxes on the 1-st column from the top. 
Then take all vectors in these boxes and form a vector subspace 
$V_1$ generated by them.  In the above example, $q_1 = 3$; hence 
$V_1 = \mathbf{C}<e_1, e_6, e_{10}>$.  We next delete these $q_1$ 
boxes from the original Young tablaux to get a new one. The 
new tablaux has exatly $d_i -1$ boxes on the $i$-th row for  
$1 \le i \le q_1$, and has exactly $d_i$ boxes on the $i$-th row for 
$i > q_1$:   

\begin{picture}(100, 100)(0, 0)
\put(00,  80){\line(1, 0){80}}
\put(00,  60){\line(1, 0){80}}
\put(00,  40){\line(1, 0){60}}
\put(00,  20){\line(1, 0){60}}
\put(00,  00){\line(1, 0){20}}
\put(00,  00){\line(0, 1){80}}
\put(20,  00){\line(0, 1){80}}
\put(40,  20){\line(0, 1){60}}
\put(60,  20){\line(0, 1){60}}
\put(80,  60){\line(0, 1){20}}
\put(05, 65){$e_2$}
\put(25, 65){$e_3$}
\put(45, 65){$e_4$}
\put(65, 65){$e_5$}
\put(05, 45){$e_7$}
\put(25, 45){$e_8$}
\put(45, 45){$e_9$}
\put(05, 25){$e_{11}$}
\put(25, 25){$e_{12}$}
\put(45, 25){$e_{13}$}
\put(05, 05){$e_{14}$}
\end{picture}

Consider the $q_2$ boxes on the 1-st column of the new diagram from 
the top and take all vectors in these boxes. They and $V_1$ together 
generate a vector subspace $V_2$.  In the above example, 
$V_2 = \mathbf{C}<e_1, e_6, e_{10}, e_2, e_7, e_{11}, e_{14}>$. 
Deleting the $q_2$ boxes, we get again a new Young tablaux. 
Repeat the similar process and we finally get a desired flag 
$V_1 \subset V_2 \subset ... \subset V_l = \mathbf{C}^n$ of 
type $(q_1, ..., q_l)$. One can check that this is a unique flag 
of type $(q_1, ..., q_l)$ such that $y(V_i) \subset V_{i-1}$. 
Therefore, $\nu^{-1}(y)$ consists of one point. 

{\bf Example (2.5.4)}.  
By $J'$ in (2.2) we introduce a non-degenerate skew symmetric form $<\;,  \;>$ on ${\bf C}^4$. 
Define $$SP(4) := \{A \in GL(4, \mathbf{C}); A^tJ'A = J'\}.$$ Its Lie algebra 
is $sp(4)$.  By an easy calculation, we see that   
$$ x = \left( \begin{array}{cccc}
0 & 1 & 0 & 0 \\ 
0 & 0 & 0 & 0 \\ 
0 & 0 & 0 & -1 \\
0 & 0 & 0 & 0 
\end{array}\right) $$ is an element of $so(4)$.  Note that $x$ has Jordan type 
$[2, 2]$.  In general, a parabolic subgroup of 
$SP(2n)$ is obtained as the group of stabilizers of an isotropic flag 
$\{V_i\}_{1 \leq i \leq s}$ of ${\mathbf C}^4$. An isotropic flag is a flag 
such that $V_i^{\perp} = V_{s-i}$ for all $i$.  
The (flag) type of an isotropic flag can be written as $(p_1,  . . . ,  p_k, q, p_k,  . . . ,  p_1)$ 
with some positive integers $p_i$ and a non-negative integer $q$.  Here we put $q = 0$ 
if the length of the flag is even.  Let $e_1$, ..., $e_4$ be the standard base of $\mathbf{C}^4$.  
Since $V := \mathbf{C}<e_1,  e_3>$ is a Lagrangian subspace (i.e. $V^{\perp} = V$),  the flag 
$V \subset \mathbf{C}^4$ is isotropic of type  (2, 2).  Let 
$Q_{2, 2}$ be the stabilizer group of this isotropic flag.  Since   
$x\cdot {\mathbf C}^4 \subset V$,  $x\cdot V = 0$,  we have $x \in \n({\q}_{2, 2})$. 
Since $\dim O_{[2,2]} = 2 \dim SP(4)/Q_{2,2}$, we know that $O_[2,2]$ is the 
Richardson orbit for $Q_{2,2}$. The Springer map 
$$ \nu_{2,2}: SP(4) \times^{Q_{2,2}} \n({\q}_{2, 2}) \to \bar{O}_{[2,2]} $$ 
is a birational map. 

On the other hand, let us consider the isotropic flag:  
$$ V_1 := \mathbf{C}<e_1>, \; V_2 := \mathbf{C}<e_1,  e_2,  e_3>.  $$ 
Let $Q_{1,2,1} \subset SP(4)$ be the stabilizer group of this flag.  
Since $\dim O_{[2,2]} = 2 \dim SO(4)/Q_{1,2,1}$ and $x \in \n(\q_{1,2,1})$, 
we see that $O_{[2,2]}$ is the Richardson orbit for $Q_{1,2,1}$. 
The Springer map 
$$ \nu_{1,2,1}: SP(4) \times^{Q_{1,2,1}} \n({\q}_{1,2, 1}) \to \bar{O}_{[2,2]} $$
is {\em not} birational. In fact, let us consider the isotropic flag  
$$ V'_1 := \mathbf{C}<e_3>, \; V'_2 := \mathbf{C}<e_1,  e_3,  e_4>, $$ and 
its stabilizer group $Q'_{1,2,1}$. Then $x \in \n(\q'_{1,2,1})$ and 
$Q_{1,2,1}$ and $Q'_{1,2,1}$ are conjugate to each other. 
This means that $\nu_{1,2,1}^{-1}(x)$ contains at least two points. One can 
prove that $\mathrm{deg}\nu_{1,2,1} = 2$ (cf. [He]).    
\vspace{0.2cm}

{\bf Example (2.5.5)}. Assume that $\g = sp(m)$ or $so(m)$. Let $z \in \g$ be a 
nilpotent element of Jordan type $\mathbf{d} := [d_1^{s_1}, d_2^{s_2}, ..., d_k^{s_k}]$. 
Let $O_z$ be the nilpotent orbit containing $z$. 
Assume that $d_p \geq d_{p+1} + 2$ for some $p$. Put $r := \Sigma_{1 \le j \le p}s_j$. 

We shall show that there are a parabolic subgroup $Q$ of $G$  
with flag type $(r, m-2r,r)$, a Levi decomposition $\q = \l \oplus \n$, and a nilpotent 
orbit $O'$ of $\l$ such that $O_z = \mathrm{Ind}^{\g}_{\l}(O')$. 
Here $$\l = gl(r) \oplus \g',$$ where $\g' = sp(m)$ (resp. $so(m)$) if  
$\g = sp(m)$ (resp. $so(m)$). The orbit $O'$ is a nilpotent orbit of $\g'$ 
with Jordan type 
$\mathbf{d}' := [(d_1-2)^{s_1}, ..., (d_p -2)^{s_p}, d_{p+1}^{s_{p+1}}, ..., d_k^{s_k}]$.
 
Let us consider the case $\g = sp(m)$.    
We prepare two skew-symmetric 
vector space $V_d$ ($d$: even), and $W_{2d}$ ($d$: odd) as follows. 
The vector space $V_d$ is a $d$-dimensional vector space with a skew-symmetric form determined by the  
$d \times d$ matrix $J'$ in (2.2).   
Let $Z_d$ be the $d \times d$ matrix such that 
$Z_d(i,i+1) = 1$ ($1 \leq i \leq d/2$), $Z_d(i,i+1) = -1$ ($d/2 + 1 \leq i \leq d-1$) 
and otherwise $Z_d(i,j) = 0$. We denote by $z_d$ the endomorphism of $V_d$ 
determined by $Z_d$. 
The vector space $W_{2d}$ is a $2d$ dimensional vector space with a skew-symmetric form determined 
by the $2d \times 2d$ matrix $J'$ in (2.2). By using the Jordan matrix $J_d$ we define   
$$ Z_{2d} := \left( \begin{array}{cccc}
J_{d} & 0 \\ 
0 & -J_{d}  
\end{array}\right)$$  
and let $z_{2d}$ be the corresponding endomorphism of $W_{2d}$.  
  
Note that, in the partition $\mathbf{d}$,  $s_i$ is even if $d_i$ is odd.  
When $d_i$ is even, we put  $U_i := V_{d_i}^{\oplus s_i}$ and 
define $z_i \in \mathrm{End}(U_i)$ by $z_i = z_{d_i}^{\oplus s_i}$. 
When $d_i$ is odd, we put $U_i := W_{2d_i}^{\oplus s_i/2}$ and 
define $z_i \in \mathrm{End}(U_i)$ by $z_i = z_{2d_i}^{\oplus s_i/2}$. 
We may assume that   
$$ (\mathbf{C}^m, <\;, \;>) := \oplus_{1 \le i \le k} U_i,$$ and 
$z = \oplus z_i$. 
Each $U_i$ has a filtration $0 \subset U_{i,1} \subset U_{i,2} 
\subset ... \subset U_{i, d_i} = U_i$ defined by $U_{i,j} := \mathrm{Im}(z_i^{d_i-j})$.  
We put $$ F := \oplus_{1 \le i \le p}U_{i,1}.$$ 
By definition $\dim F = r$ and $F \subset F^{\perp}$. 
Moreover, one can check that $z\vert_{F^{\perp}/F}$ is a nilpotent endomorphism 
of Jordan type $\mathbf{d}' = [(d_1-2)^{s_1}, ..., (d_p -2)^{s_p}, d_{p+1}^{s_{p+1}}, ..., d_k^{d_k}]$. 
Conversely, $F$ is the unique $r$ dimensional isotropic subspace such that 
$z(F) = 0$ and $z\vert_{F^{\perp}/F}$ has Jordan type $\mathbf{d}'$. 
Let $Q \subset SP(m)$ be the stabilizer group of the isotropic flag 
$$ 0 \subset F \subset F^{\perp} \subset \mathbf{C}^m,$$ and let $O'$ be the nilpotent 
orbit of $sp(F^{\perp}/F)$ which contains $z\vert_{F^{\perp}/F}$. 
Then $O_z = \mathrm{Ind}^{\g}_{\l}(O')$. Moreover, the generalized Springer map 
$$ \nu_{Q}: SP(m) \times^Q(\n + \bar{O'}) \to \bar{O}_z $$ is birational.       
The case $\g = so(m)$ is similar. 
\vspace{0.2cm}

{\bf Example (2.5.6)}. Let us consider a nilpotent element $x \in so(4n+2)$ 
with Jordan type $[2^{2n}, 1^2]$. Denote by $O_x$ the nilpotent 
orbit containing $x$. Let $V_2$ be a 2-dimesional 
vector space with a symmetric form determined by the matrix 
$$\left( \begin{array}{cccc}
0 & 1 \\ 
1 & 0  
\end{array}\right).$$ 
Let $W_4$ be a 4-dimensional vector space with a 
symmetric form determined by 
$$\left(\begin{array}{cccc} 
0 & 0 & 0 & 1 \\
0 & 0 & 1 & 0 \\ 
0 & 1 & 0 & 0 \\ 
1 & 0 & 0 & 0
\end{array}\right).$$ 
Define $z \in \mathrm{End}(W_4)$ by 
the matrix 
$$\left(\begin{array}{cccc} 
0 & 1 & 0 & 0 \\
0 & 0 & 0 & 0 \\ 
0 & 0 & 0 & -1 \\ 
0 & 0 & 0 & 0
\end{array}\right).$$
One may assume that $$ x = z^{\oplus n} \oplus 0 \in  W_4^{\oplus n} \oplus V_2.$$ 
Let $e^{(i)}_1$, ..., $e^{(i)}_4$ be the (standard) basis of the $i$-th direct summand $W_4$ 
of $W_4^{\oplus n}$, and let $f_1$ and $f_2$ be the basis of $V_2$. 
Define a $2n+1$-dimensional isotropic subspace $F$ by 
$$ F := \mathbf{C}<e^{(1)}_1, e^{(1)}_3, e^{(2)}_1, e^{(2)}_3, ..., e^{(n)}_1, e^{(n)}_3, f_1> $$ 
and consider the istropic flag $\{F_{\cdot}\}$ defined by $F_0 = 0$, $F_1 := F$, $F_2 := F^{\perp}$ 
and $F_3 = \mathbf{C}^{4n+2}$. One can check that $x(F_i) \subset F_{i-1}$ for all $i$. 
Let $Q \subset SO(4n+2)$ be the parabolic subgroup stabilizing this flag. 
One can check that $x \in \n(\q)$ and $\dim O_x = 2 \dim SO(4n+2)/Q$. 
Therefore, $O_x$ is the Richardson orbit for $Q$.  
But $\{F_{\cdot}\}$ is not the unique isotropic flag (of type $(2n+1, 2n+1)$) with 
this property.  We put  
$$ F' := \mathbf{C}<e^{(1)}_1, e^{(1)}_3, e^{(2)}_1, e^{(2)}_3, ..., e^{(n)}_1, e^{(n)}_3, f_2>, $$
and define $\{F'_{\cdot}\}$ by $F'_0 = 0$, $F'_1 := F'$, $F'_2 := (F')^{\perp}$ 
and $F'_3 := \mathbf{C}^{4n+2}$. Then $\{F'_{\cdot}\}$ has the same property. 
Let $Q'$ be the corresponding parabolic subgroup of $SO(4n+2)$. Then $O_x$ 
is the Richardson orbit for $Q'$. Although $Q$ and $Q'$ have the same flag 
types, $Q$ and $Q'$ are {\em not conjugate} as the subgroups of $SO(4n+2)$. 
As a consequence, we have two Springer map 
$$ SO(4n+2) \times^Q \n(\q) \to \bar{O}_x \leftarrow SO(4n+2) \times^{Q'} \n(\q') $$ 
and both of them are birational. @

(2.6) ({\em Crepant resolutions and {\bf Q}-factorial terminalizations}): 
As before, let $O$ be a nilpotent orbit of $\g$ and let $\tilde{O}$ be the 
normalization of $\bar{O}$. It is an important problem to find a crepant 
resolution or its substitute for $\tilde{O}$. Let us recall 

{\bf Definition (2.6.1)}. Let $X$ and $Y$ be normal variety with rational Gorenstein 
singularities and let $\mu: Y \to X$ be a birational projective morphism. 
Then $\mu$ is called a {\em crepant resolution} (resp. {\em $\mathbf{Q}$-factorial 
terminalization}) of $X$ if $Y$ is smooth (resp. $Y$ has only {\bf Q}-factorial 
terminal singularities) and $K_Y = \mu^*K_X$. 
\vspace{0.2cm}

{\bf Theorem (2.6.2)}. {\em Let $O$ be a nilpotent orbit of a complex simple Lie 
algebra $\g$. Then there are a parabolic subalgebra $\q$ of $\g$ and a nilpotent orbit $O'$ of 
$\l(\q)$ such that the following holds}: 
\vspace{0.15cm}
 
(1) $O = \mathrm{Ind}^{\g}_{\l}(O')$. 

(2) {\em Let $\nu: G \times^Q (\n + \bar{O}') \to \bar{O}$ be 
the generalized Springer map. Then its normalized map $\nu^n$ (cf. (2.5.1)) 
gives a {\bf Q}-factorial terminalization of $\tilde{O}$.} 
\vspace{0.15cm}
 
Theorem (2.6.2) is due to [Na 3] when $\g$ is of 
classical type and is due to [Fu 2] when $\g$ is of exceptional 
type. Here we shall give a rough sketch of the proof when $\g$ is 
classical. First let us consider the case $\g = sl(n)$. Assume that $O$ has Jordan type 
$\mathbf{d}$. Let $\mathbf{d}^t = [q_1, q_2, ..., q_l]$ be the dual partition 
of $\mathbf{d}$. By Example (2.5.1), $O$ is the Richardson orbit for a parabolic 
subgroup $Q \subset SL(n)$ of flag type $(q_1, q_2, ..., q_l)$. Moreover, the Springer 
map $\nu: G \times^Q \n \to \bar{O}$ is birational. Note that $G \times^Q \n$ is isomorphic 
to the cotangent bundle $T^*(G/Q)$ of the homogeneous space $G/Q$. The pull-back of 
the Kostant-Kirillov 2-form $\omega$ (cf. (2.1)) coincides with the canonical 2-form 
of $T^*(G/Q)$; hence $\nu$ is a crepant resolution. 
Next let us consider the cases $\g = sp(m)$ and $\g = so(m)$. 
We say that a partition $\mathbf{d} := [d_1^{s_1}, d_2^{s_2}, ..., d_k^{s_k}]$ of $m$ 
has full members if $d_i = k + 1 -i$ for all $i$. 

{\bf Proposition (2.6.3)}. {\em Assume that $\g = sp(m)$ or $so(m)$. Then 
$\tilde{O}_{\mathbf{d}}$ has terminal singularities if and only if $\mathbf{d}$ 
has full members. If $\mathbf{d}$ has full members, then $\tilde{O}_{\mathbf d}$ 
is {\bf Q}-factorial except when $\g = so(4n+2)$, $n \geq 1$ and $\mathbf{d} = 
[2^{2n}, 1^2]$.} 
\vspace{0.2cm}

For the proof of Proposition (2.6.3), see [Na 3]. 
\vspace{0.15cm}

Assume that $\tilde{O}$ does not have {\bf Q}-factorial terminal singularities. 
By (2.6.3), $O$ does not have full members or $O = O_{[2^{2n}, 1^2]} \subset so(4n+2)$. 
In the second case, $O$ is a Richardson orbit and has a crepant resolution by 
(2.5.6). In the first case, let $\mathbf{d} := [d_1^{s_1}, ..., d_k^{s_k}]$ be the Jordan 
type of $O$. Then $d_p \geq d_{p+1} + 2$ for some $p$. 
The situation is now the same as (2.5.5).  The orbit  
$O$ is induced, and as in (2.5.5) one can find a generalized 
Springer map $$\nu: G \times^Q (\n(\q) + \bar{O'}) \to \bar{O}, $$ which is 
birational. Then $O'$ is again a nilpotent orbit of a smaller Lie algebra of the same type. 
If $O'$ already has {\bf Q}-factorial terminal singularities, then 
the normalization of  $G \times^Q (\n(\q) + \bar{O'})$ gives a {\bf Q}-factorial 
terminalization (cf. Proposition (4.2) below,  see also [Na 3, Lemma (1.2.4)]
\footnote{The proof of [Na 3], Lemma (1.2.4) contains an error. 
In fact, $<y, [v_1, w_1]>$ is claimed to be zero there, but it is not correct. 
The equality starting from line 9, p 552 should contain the additional term 
$<y, [v_1,w_1]>$. The equality on line -2, p 552 should also contain 
$<y, [v_1,w_1]>$. But the claim itself is correct.}). 
If not, then we repeat the same process. 
By (2.5.2) we have a 
birational map $$\nu': G \times^{Q'} (\n(\q') + \bar{O}'') \to 
G \times^Q (\n(\q) + \bar{O'}).$$ Finally we get a {\bf Q}-factorial terminalization 
of $\tilde{O}$. 
\vspace{0.2cm}

{\bf  \S 3. Birational geometry of Q-factorial terminalizations}
 
(3.1) {\em Parabolic subgroups and root systems}:  
Let $G$ be a simple algebraic group over {\bf C} and let $\g$ be its Lie algebra.   
We fix a maximal torus $T$ of $G$ and denote by 
$\h$ its Lie algebra.  Let $\Phi$ be the root system for $\g$ determined by 
$\h$ (cf. [Hu 1]).  The root system $\Phi$ has a natural involution $-1$.  
There is a (unique) involution $\varphi_{\g}$ of $\g$ which stabilizes 
$\h$ and which acts on $\Phi$ via $-1$ (cf. [Hu 1], 14.3).  
Let 
$$\g = \h \oplus \bigoplus_{\alpha \in \Phi}\g_{\alpha}$$
be the root space decomposition. 
Let us choose a base $\Delta$ of $\Phi$ and denote by 
$\Phi^+$ the set of positive roots. Then 
$\b = \h \oplus 
\bigoplus_{\alpha \in \Phi^+}\g_{\alpha}$ is a Borel subalgebra 
of $\g$ which contains $\h$.  Let $B$ be the corresponding 
Borel subgroup of $G$. 
Take a subset $I$ of $\Delta$.  Let $\Phi_I$ be the root 
subsystem of $\Phi$ generated by $I$ and put 
$\Phi^-_I := \Phi_I \cap \Phi^-$, where $\Phi^{-}$ is the 
set of negative roots. Then 
$$\q_I := \h \oplus\bigoplus_{\alpha \in \Phi^-_I}\g_{\alpha} 
\oplus \bigoplus_{\alpha \in \Phi^+}\g_{\alpha}$$ 
is a parabolic subalagebra containing $\b$. 
This parabolic subalgebra $\q_I$ is called a standard parabolic 
subalgebra with respect to $I$.  Let $Q_I$ be the corresponding 
parabolic subgroup of $G$. By definition, $B \subset Q_I$. 
Any parabolic subgroup $Q$ of $G$ is conjugate to a 
standard parabolic subgroup $Q_I$ for some $I$.  Moreover, if 
two subsets $I$, $I'$ of $\Delta$ are different, $Q_I$ and 
$Q_{I'}$ are not conjugate. Thus, a conjugacy class of parabolic 
subgroups of $G$ is completely determined by $I \subset \Delta$. 
In this paper, to specify the subset $I$ of $\Delta$, we shall 
use the {\em marked} Dynkin diagram. Recall that $\Delta \subset \Phi$ 
defines a Dynkin diagram; each vertex corresponds to a simple root (an 
element of $\Delta$).  Now, if a subset $I$ of $\Delta$ is given, we 
indicate the vertices corresponding to $I$ by white vertices, and 
other vertices by black vertices.  A black vertex is called a {\em marked 
vertex}.  A Dynkin diagram with such a marking is called a {\em 
marked Dynkin diagram}, and a marked Dynkin diagram with only one 
marked vertex is called a {\em single} marked Dynkin diagram. 
Note that the standard parabolic subgroup corresponding to a 
single marked Dynkin diagram (resp. full marked Dynkin diagram)  
is a maximal parabolic subgroup (resp. a Borel subgroup).  
Let $\q$ be a parabolic subalgebra of $\g$ which contains $\h$. 
Let $\r(\q)$ (resp. $\n(\p)$) be the solvable radical (resp. nilpotent 
radical) of $\q$.  
We put $\k(\q) := \r(\q) \cap \h$. Then $$\r(\q) = \k(\q) \oplus \n(\q).$$  
On the other hand, the Levi factor $\l(\q)$ of $\q$ is defined 
as $\l(\q) := \g^{\k(\q)}$. Here, 
$\g^{\k(\q)} := \{x \in \g; [x, y] = 0, \forall y \in \k(\q)\}. $ 
Note that $\k(\q)$ is the center of $\l(\q)$.  
Then $$\q = \l(\q) \oplus \n(\q)$$ and 
$$\l(\q) = \k(\q) 
\oplus [\l(\q), \l(\q)].$$ 
If $\q = \q_I$, we have 
$$ \k(\q_I) = \{h \in \h; \alpha(h) = 0, \forall \alpha \in I\}.$$ 
Moreover, we define 
$$ \k(\q_I)^{reg} = \{h \in \k(\q_I); \alpha (h) \ne 0, \forall \alpha \in \Phi\setminus 
\Phi_I\}. $$ Note that $\k(\q_I)^{reg}$ is an open subset of 
$\k(\q_I)$. 

(3.2) ({\em Parabolic subalgebras with a fixed Levi part}):  
A subalgebra $\l$ of $\g$ is called a Levi subalgebra if it is the 
Levi part of some parabolic subalgebra $\q$ of $\g$. 
We fix a Levi subalgebra $\l$ and put 
$$ \mathcal{S}(\l) := 
\{\mathrm{parabolic}\;  \mathrm{subalgebras}\;  \q \; \mathrm{of} \; \g ;\; \l(\q) = \l \}. $$
Fix a Cartan subalgebra $\h$ of $\g$ so that $\h \subset \l$. 
Let $\b$ be a Borel subalgebra of $\g$ so that $\h \subset \b \subset \q$. 
Then $\q$ corresponds to a 
marked Dynkin diagram $D$. Take a marked vertex $v$ of 
the Dynkin diagram $D$ and consider the maximal connected 
single marked Dynkin subdiagram $D_v$ of $D$ containing $v$. 
We call $D_v$ the {\em single marked diagram associated with 
$v$}.  When $D_v$ is one of the following, 
we say that $D_v$ (or $v$) {\em is of the 
first kind}, and when $D_v$ does not coincide with any of them, 
we say that $D_v$ (or $v$) {\em is of the second kind}. 
\vspace{0.6cm}

$A_{n-1}$ $(k < n/2)$ 
  
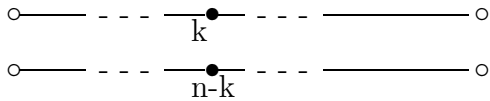
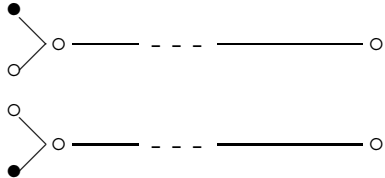
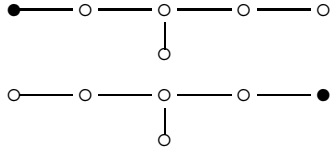
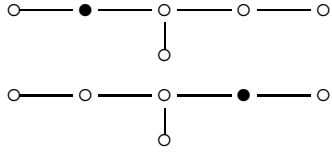
\begin{picture}(300,20)(0,0) 
\put(30,-3){$\circ$}\put(35,0){\line(1,0){25}} 
\put(65,-3.5){- - -}\put(90,0){\line(1,0){15}}
\put(105,-3){$\bullet$}\put(110,0){\line(1,0){10}}
\put(100,-10){k}\put(125,-3.5){- - -}\put(150,0)
{\line(1,0){55}}\put(207,-3){$\circ$}    
\end{picture}  

\begin{picture}(300,20)(0,0) 
\put(30,-3){$\circ$}\put(35,0){\line(1,0){25}} 
\put(65,-3.5){- - -}\put(90,0){\line(1,0){15}}
\put(105,-3){$\bullet$}\put(110,0){\line(1,0){10}}
\put(100,-10){n-k}\put(125,-3.5){- - -}\put(150,0)
{\line(1,0){55}}\put(207,-3){$\circ$}    
\end{picture}
\vspace{0.4cm}

$D_n$ $(n:$ $\mathrm{odd} \geq 5)$  

\begin{picture}(300,20)(0,0) 
\put(30,10){$\bullet$}\put(30,-13){$\circ$}
\put(35,10){\line(1,-1){10}}\put(35,-10){\line(1,1){10}}
\put(47,-3){$\circ$}\put(55,0){\line(1,0){25}}
\put(85,-3.5){- - -}\put(110,0){\line(1,0){55}}\put(167,-3)
{$\circ$}
\end{picture} 
\vspace{0.6cm}

\begin{picture}(300,20)(0,0) 
\put(30,10){$\circ$}\put(30,-13){$\bullet$}
\put(35,10){\line(1,-1){10}}\put(35,-10){\line(1,1){10}}
\put(47,-3){$\circ$}\put(55,0){\line(1,0){25}}
\put(85,-3.5){- - -}\put(110,0){\line(1,0){55}}\put(167,-3)
{$\circ$}
\end{picture}       
\vspace{0.6cm}

$E_{6,I}$:  

\begin{picture}(300,20)
\put(30,-3){$\bullet$}\put(35,0){\line(1,0){20}}
\put(57,-3){$\circ$}\put(65,0){\line(1,0){20}}
\put(87,-3){$\circ$}\put(90,-5){\line(0,-1){10}}
\put(87,-20){$\circ$}\put(95,0){\line(1,0){20}}
\put(117,-3){$\circ$}\put(125,0){\line(1,0){20}}
\put(147,-3){$\circ$} 
\end{picture} 
\vspace{0.4cm}

\begin{picture}(300,20)
\put(30,-3){$\circ$}\put(35,0){\line(1,0){20}}
\put(57,-3){$\circ$}\put(65,0){\line(1,0){20}}
\put(87,-3){$\circ$}\put(90,-5){\line(0,-1){10}}
\put(87,-20){$\circ$}\put(95,0){\line(1,0){20}}
\put(117,-3){$\circ$}\put(125,0){\line(1,0){20}}
\put(147,-3){$\bullet$} 
\end{picture}
\vspace{0.4cm}

$E_{6,II}$: 

\begin{picture}(300,20)
\put(30,-3){$\circ$}\put(35,0){\line(1,0){20}}
\put(57,-3){$\bullet$}\put(65,0){\line(1,0){20}}
\put(87,-3){$\circ$}\put(90,-5){\line(0,-1){10}}
\put(87,-20){$\circ$}\put(95,0){\line(1,0){20}}
\put(117,-3){$\circ$}\put(125,0){\line(1,0){20}}
\put(147,-3){$\circ$} 
\end{picture}
\vspace{0.4cm}

\begin{picture}(300,20)
\put(30,-3){$\circ$}\put(35,0){\line(1,0){20}}
\put(57,-3){$\circ$}\put(65,0){\line(1,0){20}}
\put(87,-3){$\circ$}\put(90,-5){\line(0,-1){10}}
\put(87,-20){$\circ$}\put(95,0){\line(1,0){20}}
\put(117,-3){$\bullet$}\put(125,0){\line(1,0){20}}
\put(147,-3){$\circ$} 
\end{picture}

\vspace{1.0cm}

In the single marked Dynkin diagrams above, two diagrams 
in each type (i.e. $A_{n-1}$, $D_n$, $E_{6,I}$, $E_{6,II}$) 
are called {\em duals}. The Weyl group 
$W$ of $\g$ does not contain $-1$ exactly when $\g = A_n (n \geq 2)$, $D_n$ (n: odd) or 
$E_6$ (cf. [Hu 1], p.71, Exercise 5). This property characterizes the  
Dynkin diagrams in the list.  Moreover, the single marked Dynkin diagrams $D_v$ in the list are  
characterized by the following property. 
\vspace{0.15cm}

(*) Let $\q_v$ be the parabolic       
subalgebra of $\g$ corresponding to $D_v$, and let $\varphi_{\g}$ be the automorphism of 
$\g$ determined by $-1$ (cf.(3.1)). 
Then $\varphi_{\g}(\q_v)$ is not   
conjugate to $\q_v$. If $\varphi_{\g}(\q_v)$ corresponds to a single marked Dynkin 
diagram $D^*_v$, then $D_v$ and $D^*_v$ are mutually duals. 
\vspace{0.15cm}
 
Let $\bar{D}$ be the marked Dynkin diagram obtained from 
$D$ by making $v$ unmarked. 
Let $\bar{\q}$ be the parabolic subalgebra containing 
$\q$ corresponding to $\bar{D}$. Now let us define a new 
marked Dynkin diagram $D'$ as follows. If $D_v$ is of the 
first kind, we replace $D_v \subset D$ by its dual diagram 
$D_v^*$ to get a new marked Dynkin diagram $D'$. If $D_v$ is 
of the second kind, we define $D' := D$. 
The new diagram $D'$ obtained in this way is called an {\em 
adjacent diagram} to $D$. 
As in (3.1),  the set of unmarked vertices 
of $D$ (resp. $\bar{D}$) defines a subset 
$I \subset \Delta$ (resp. $\bar{I} \subset 
\Delta$). 
By definition, $v \in \bar{I}$. The unmarked vertices of 
$\bar{D}$ define a Dynkin subdiagram, which is decomposed 
into the disjoint sum of the connected component containing $v$ and 
the union of other components. Correspondingly, we have a decomposition 
$\bar{I} = I_v \cup I'_v$ with $v \in I_v$.  
The parabolic subalgebra $\q$ (resp. $\bar{\q}$) coincides 
with the standard parabolic subalgebra $\q_I$ (resp. 
${\q}_{\bar{I}}$). Let $\mathfrak{l}_{\bar{I}}$ be the (standard) 
Levi factor of $\q_{\bar{I}}$. Let $\z(\l_{\bar{I}})$ be the center 
of $\l_{\bar{I}}$. Then $\l_{\bar{I}}/\z(l_{\bar{I}})$ is decomposed 
into the direct sum of simple factors. Now let $\l_{I_v}$ be the 
simple factor corresponding to $I_v$ and let $\l_{I'_v}$ be the 
direct sum of other simple factors. Then 
$$ \mathfrak{l}_{\bar{I}}/\z(\l_{\bar{I}}) = \mathfrak{l}_{I_v} \oplus 
\mathfrak{l}_{I'_v}. $$ 
The marked Dynkin diagram $D_v$ defines a standard parabolic subalgebra 
$\q_v$ of $\mathfrak{l}_{I_v}$. Here let us consider the involution  
$\varphi_{\mathfrak{l}_{I_v}} \in \mathrm{Aut}(\mathfrak{l}_{I_v})$ (cf. (3.1)).      
When $D_v$ is of the first kind, $\varphi_{\mathfrak{l}_{I_v}}(\q_v)$ is 
conjugate to a standard parabolic subalagebra of 
$\mathfrak{l}_{I_v}$ with the dual marked Dynkin diagram 
$D_v^*$ of $D_v$. When $D_v$ is of the second kind, $\varphi_{\mathfrak{l}_{I_v}}(\q_v)$ 
is conjugate to $\q_v$ in $\mathfrak{l}_{I_v}$.  
Let $\pi: \l_{\bar{I}} \to \l_{\bar{I}}/\z(\l_{\bar{I}})$ be the quotient homomorphism.  
Note that 
$$ \bar{\q} = \mathfrak{l}_{\bar{I}} \oplus n(\bar{\q}), $$ 
$$ \q = {\pi}^{-1}(\q_v \oplus \mathfrak{l}_{I'_v}) \oplus n(\bar{\q}). $$ 
Here we define 
$$ \q' = {\pi}^{-1}(\varphi_{\mathfrak{l}_{I_v}}(\q_v) \oplus \mathfrak{l}_{I'_v}) \oplus 
n(\bar{\q}). $$ 
Then $\q' \in \mathcal{S}(\l)$ and $\q'$ is conjugate 
to a standard parabolic subalgebra with the marked Dynkin 
diagram $D'$. {\em This $\p'$ is said to 
be the parabolic subalgebra twisted by $v$}. 
Two marked diagrams $D_1$ and $D_2$ are called {\em equivalent} if there is a 
finite chain of adjacent diagrams connecting $D_1$ and $D_2$. 
\vspace{0.2cm}

{\bf Definition (3.2.1)} ({\em Mukai flops and primitive pairs}): 
Let $\q$ and $\q'$ be two parabolic subalgebras of $\g$ corresponding to 
the dual diagrams in the list above.  Assume that $\q$ and $\q'$ have a common Levi 
factor $\l$. When the diagram is of type $A_{n-1,k}$, $D_n$ or 
$E_{6,II}$, define $O'$ to be the 0-orbit in $\l$. 
When the diagram is of type $E_{6,I}$, define $O'$ to be the 0-orbit, 
$O_{[3,2^2,1^3]}$ or $O_{[2^2,1^6]}$. Such a pair $(\q, O')$ is called a 
{\em primitive pair}. 
We put $O := \mathrm{Ind}^{\g}_{\l}(O')$. Then we have a  
diagram of normalized maps (cf. (2.5.1)) of  (generalized) Springer maps:  
$$ G \times^Q (\n(\q) + \tilde{O'}) \stackrel{\nu^n}\to \tilde{O} \stackrel{(\nu')^n}\leftarrow 
G \times^{Q'} (\n(\q') + \tilde{O'}).$$ For every primitive pair, $\nu^n$ and $(\nu')^n$ 
are both isomorphisms in codimension one. 
This diagram is called a {\em Mukai flop}.
\vspace{0.2cm}

(3.3) ({\em Brieskorn-Slodowy diagram}): Let $O \subset \g$ be the induced orbit from $O' \subset \l(\q)$ for a 
parabolic subalgebra $\q$.  Assume that $\tilde{O'}$ has only {\bf Q}-factorial terminal singularities.     
Let us consider the subvariety $$\r(\q) + \bar{O}' \subset \r(\q) \oplus [\l(\q), \l(\q)].$$ 
The $Q$-adjoint action stabilizes $\r(\q) + \bar{O}'$ as a set; hence one has an associated 
fiber bundle $$X'_{\q, O'} := G \times^Q (\r(\q) + \bar{O}')$$ over $G/Q$. 
Write $x \in \r(\q) + \bar{O}'$ as $x = x_1 + x_2 + x_3$, where $x_1 \in \k(\q)$, 
$x_2 \in \n(\q)$ and $x_3 \in [\l(\q), \l(\q)]$. Define a map 
$$\eta: X'_{\q, O'}  \to \k(\q) $$ by $[g,x] \in X_{\q} \to x_1$. This is a well-defined map; in fact, 
for $q \in Q$, we have $(Ad_q(x_1))_1 = x_1$, $Ad_q(x_2) \in \n(\q)$ and $Ad_q(x_3) \in 
\n(\q) \oplus [\l(q), \l(\q)]$. 
\vspace{0.2cm}

{\bf Lemma (3.3.1)}. {\em For $t \in \k(\q)^{reg}$, any orbit of 
the $Q$-variety $t + \n(q) + \bar{O}'$ is of the form $Q(t + y)$ 
with $y \in \bar{O'}$.}   

{\em Proof}. We shall prove that $\cup_{y \in \bar{O}'}Q\cdot (t+y) = t + \n(\q) + \bar{O}'$.  
Define $Z_Q(t+y) := \{q \in Q; Ad_q(t+y) = t+y\}$. Then 
$Q\cdot (t+y) \cong Q/Z_Q(t+y)$. 
Note that $t$ (resp. $y$) is the 
semi-simple part (resp. nilpotent part) of $t+y$ in the Jordan-Chevalley 
decomposition because $[t,y] = 0$. Hence $Z_Q(t+y) = Z_Q(t) \cap Z_Q(y)$. Since $t \in \k(\q)^{reg}$, 
we have $Z_Q(t) = L(Q)$, and $Z_Q(t) \cap Z_Q(y) = Z_{L(Q)}(y)$. Let $O_y \subset \l(\q)$ 
be the $L(Q)$-adjoint orbit containing $y$.  Then one has  
$$\dim Q/Z_Q(t+y)  = \dim \n(\q) + \dim O_y.$$ Let us write $Q = U(Q)\cdot L(Q)$ 
with the unipotent radical $U(Q)$. 
Since $L(Q)\cdot (t+y) = t + O_y$ and $U(Q)\cdot (t+y) \subset t + y + \n(\q)$, we 
see that $U(Q)\cdot (t+y)$ is dense in $t + y + \n(\q)$.  
But any $U(Q)$-orbit in $t + y + \n(\q)$ is closed (cf. [Hu 2], \S 17, Exercise 8), 
$U(Q)\cdot (t+y) = t + y + \n(\q)$. Q.E.D.  
\vspace{0.2cm}

We shall prove that the following diagram commutes: 
 
\begin{equation} 
\begin{CD} 
X'_{\q, O'}  @>>> G\cdot(\r(\q) + \bar{O}') \\ 
@V{\eta}VV @V{\chi}VV \\ 
\k(\q)  @>{\iota}>>  \h/W.    
\end{CD} 
\end{equation} 

Here $\chi$ is the composite of the inclusion map 
$G(\r(\q) + \bar{O}') \to \g$ and the adjoint quotient map 
$\g \to \h/W$. The horizontal map on the first row is given by 
$[g,x] \to Ad_g(x)$ and the horizontal map $\iota$ on the second row is 
the composite of two maps $\k(\q) \to \h$ and $\h \to \h/W$. 
Define $$W' := N_W(L(Q))/W(L(Q)), $$ 
where $N_W(L(Q))$ is the normalizer subgroup of $W$ for $L(Q)$. 
Then $W'$ acts on $\k(\q)$ and the normalization of $\mathrm{Im}(\iota)$ is 
$\k(\q)/W'$.  
Let us check the commutativity of the diagram. 
Choose $t \in \k(\q)^{reg}$ and $y \in O'$. Then $G \times^Q Q\cdot(t+y)$ 
is an open dense subset of $\eta^{-1}(t) = G \times^Q (t + \n(\q) + \bar{O}')$ by 
Lemma (3.3.1).   
We only have to check the commutativity of the diagram for an element $[g, t+y] \in  
G \times^Q Q\cdot(t+y)$: 

\begin{equation} 
\begin{CD} 
[g, t+y] @>>> Ad_g(t+y) \\ 
@VVV @VVV \\ 
t  @>>>  [t].    
\end{CD} 
\end{equation} 
  
The commutativity now follows from the fact that 
$Ad_g(t+y) = Ad_g(t) + Ad_g(y)$ coincides with the Jordan-Chevalley 
decomposition of $Ad_g(t+y)$. We put 
$$ Y'_{\l(\q), O'} := (\k(\q) \times_{\h/W} G\cdot(\r(\q) + \bar{O}'))_{red}. $$ 
Note that $\k(\q)$ only depends on the Levi part $\l(\q)$ of 
$\q$. Moreover, since $$G\cdot(\r(\q) + \bar{O}') = 
\overline{G\cdot (\k(\q)^{reg} + \bar{O'})},$$ $Y'_{\l(\q), O'}$ 
only depends on $\l(\q)$ and $O'$ as the index indicates.  
The commutative diagram induces a map 
$$\mu'_{\q}: X'_{\q, O'}  \to Y'_{\l(\q), O'}.$$
A nilpotent orbit of $\l(\q)$ is not necessarily stable under 
the $W'$-action. But, in our case, we have:   
\vspace{0.2cm}

{\bf Lemma (3.3.2)}: All elements $w \in W'$ stabilizes $O'$. 

{\em Proof}.  If $O' = 0$, then the statement is obvious. 
Assume that $O' \neq 0$.  Let us consider the decomposition of 
$\l(\q)$ into simple factors (up to centers). Then 
$O'$ is contained in a simple factor $\g'$ whose type in not $A$.  
Then each element $w \in W'$ induces an automorphism of Lie algebra $\g'$. 
If $\mathrm{Aut}(\g') = \mathrm{Aut}(\g')^0$, then $w$ acts on $\g$ as an 
adjoint action $Ad_g$ for some $g \in G'$. In this case, $O'$ is stable by $W'$. 
So we may assume that $\g$ is of type $D$ or $E_6$. 
Suppose that $O'$ is sent to a different nilpotent orbit $O'' \subset \g'$ by some 
$\phi \in \mathrm{Aut}(\g')$. Then $\phi$ acts on the Dynkin diagram as a 
graph automorphism.  The weighted Dynkin diagram of $O'$ should be sent 
to that of $O''$ by this graph automorphism. If $\g'$ is of type $D$, 
such things happen only  when the orbit is very even or it is $O_{[5, 1^3]} \subset 
so(8)$. Our $O'$ does not coincide with any of them.  If $\g'$ is of type 
$E_6$, then one can check that there are no such orbits by using the list of [C-M], page 
129.  Q.E.D.    
\vspace{0.2cm} 

{\bf Proposition (3.3.3)}. {\em The map $\mu'_{\q}$ is a birational projective 
morphism. In particular, $Y'_{\l(\q), O'}$ is 
irreducible. Moreover, for $t \in \k(\q)^{reg}$, the induced map 
$\eta^{-1}(t) \to \{t\} \times_{[t]}  \chi^{-1}([t])$ is a bijection.} 
\vspace{0.2cm}

{\em Proof}.  The map $\mu'_{\q}$ is written as the composite of a closed 
immersion and a projective morphism: $G \times^Q (\r(\q) + \bar{O}') 
\to G/Q \times \g \to \g$. Hence $\mu'_{\q}$ is a projective 
morphism. By (3.3.1), for $t \in \k(\q)^{reg}$, the fiber $\chi^{-1}([t])$ 
coincides with 
$$\cup_{y \in \bar{O'}} \cup_{s \in \k(\q)\mathrm{,}\;  [s] = [t]} G\cdot (s + y).$$  
But, if $[s] = [t]$, then $s = Ad_w(t)$ with some $w \in W'$.  
By (3.3.2), $Ad_w(y) \in \bar{O'}$. Therefore 
$$ \chi^{-1}([t]) =  \cup_{y \in \bar{O'}} G\cdot (t + y).$$
In particular, $\chi^{-1}([t])$ is irreducible for $t \in \k(\q)^{reg}$. 
By the argument above, any point of $\chi^{-1}([t])$ is $G$-conjugate to 
$t +  y$ with $y \in \bar{O}'$. Since $\mu'_{\q}$ is $G$-equivariant, 
it is sufficient to prove that ${\mu'}_{\q}^{-1}(t, t+y)$ consists of exactly one point,  
where $(t, t+y) \in \{t\} \times_{[t]}  \chi^{-1}([t])$. 
Assume that $[g, x] \in G \times^Q (\r(\q) + \bar{O'})$ is contained in ${\mu'}_{\q}^{-1}(t, t+y)$.  
Then $Ad_g(x) = t + y$. Moreover, by (3.1.1), one can write 
$x = Ad_q(t + y')$ with some $q \in Q$ and some $y' \in \bar{O}'$. 
This means that $Ad_{gq}(t+y') = t+y$. Since $t + y$ and $t + y'$ are 
both Jordan-Chevalley decompositions, we see that $t = Ad_{gq}(t)$ 
and $y = Ad_{gq}(y')$.  By the first equality, we have $gq \in L(Q)$, and 
hence $g \in Q$. By the second equality, we have 
$$x = Ad_q(t + y') = Ad_q(t + Ad_{q^{-1}g^{-1}}(y)) 
= Ad_q(Ad_{q^{-1}g^{-1}}(t + y)) = Ad_{g^{-1}}(t+y). $$  
As a consequence, $$[g,x] = [g, Ad_{g^{-1}}(t+y)] = [1, t+y].$$   
The rest of the argument is the same as [Na 2], Lemma 1.1. Q.E.D.          
\vspace{0.2cm}

Let $X_{\q, O'}$ be the normalization of $X'_{\q, O'} $ and let $Y_{\l(\q), O'}$ 
be the normalization of  $Y'_{\l(\q), O'}$. Then $\mu'_{\q}$ induces  
a commutative diagram 
\begin{equation} 
\begin{CD} 
X_{\q, O'} @>{\mu_{\q}}>> Y_{\l(\q), O'} \\ 
@VVV @VVV \\ 
\k(\q) @>>> \k(\q)    
\end{CD} 
\end{equation} 
Let $X_{\q,O',0}$ (resp. $Y_{\l(\q),O',0}$) be the fiber of the map $X_{\q, O'} \to \k(\q)$ 
(resp. $Y_{\l(\q),O'} \to \k(\q)$) over  
$0 \in \k(\q)$.   
\vspace{0.2cm}

{\bf Lemma (3.3.4)}. {\em Assume that $$\nu^n: G \times^Q (\n(\q) + \tilde{O'}) \to 
\tilde{O}$$ is birational
(cf. (2.5.1). Then 
one has $$X_{\q,O',0} = G \times^Q(\n(\q) + \tilde{O'}),$$  
$$Y_{\l(\q),O',0} = \tilde{O}.$$} 
{\em The map $$\mu_{\q,0}: X_{\q,O',0} \to Y_{\l(\q),O',0}$$ coincides 
with the normalized map $\nu^n$ of the generalized Springer map.} 
\vspace{0.2cm}

{\em Proof}. The first statement is obvious. 
Since $Y_{\q,O'}$ is Cohen-Macaulay and $\k(\q)$ is smooth, $Y_{\q,O',0}$ 
is also Cohen-Macaulay. The map $\mu_{\q,0}: X_{\q,O',0} \to Y_{\q,O',0}^{red}$ is 
a birational morphism with connected fibers. It factorizes the generalized 
Springer map $X_{\q,O',0} = G \times^Q (\n(\q) + \bar{O'}) \to \bar{O}$. 
By the assumption, the generalized Springer map is an isomorphism over 
$O$. This means that $\mu_{\q,0}$ is an isomorphism outside a certain 
codimension $2$ subset $Z$ of $Y_{\q,O',0}^{red}$ and $Y_{\q,O',0}^{red} - Z$ 
is smooth. Take a point $x \in 
Y_{\q,O',0}^{red} - Z$. Then we have a surjection 
$$ \mathcal{O}_{Y_{\q,O',0}} \to \mathcal{O}_{Y_{\q,O',0}^{red}} \cong 
\mathcal{O}_{X_{\q,O',0}, \mu_{\q}^{-1}(x)}.$$ 
By Nakayama's lemma, this implies that $\mathcal{O}_{Y_{\q,O'},x} 
\cong \mathcal{O}_{X_{\q,O'},\mu_{\q}^{-1}(x)}$. Therefore, $Y_{\q,O',0}$ is reduced at 
$x$, and moreover, $Y_{\q,O',0}$ is smooth at $x$. 
Since $Y_{\q,O',0}$ is Cohen-Macaulay and regular in codimension one, 
$Y_{\q,O',0}$ is normal. This means that $Y_{\q,O',0} = \tilde{O}$.   
\vspace{0.2cm}

{\bf Proposition (3.3.5)}. {\em The map $\mu_{\q}$ is crepant and is an isomorphism in 
codimension one.} 

{\em Proof}. Since $X_{\q,O',0}$ has only terminal singularities and its canonical line 
bundle is trivial, $K_{X_{\q,O'}}$ is $\mu_{\q}$-numerically trivial. By Kawamata-Viehweg 
vanishing theorem, $R^j(\mu_{\q})_*\mathcal{O}_{X_{\q,O'}} = 0$ for $j > 0$. 
Therefore, $K_{X_{\q,O'}}$ is the pull-back of a line bundle $M$ on $Y_{\l(\q),O'}$. 
Since $Y_{\l(\q),O'}$ has a $\mathbf{C}^*$-action with positive weights, its Picard 
group is trivial:  $\mathrm{Pic}(Y_{\l(\q),O'}) = 0$. This means that $K_{X_{\q,O'}}$ is a trivial 
line bundle. Then $K_{Y_{\l(\q),O'}} = (\mu_{\q})_*K_{X_{\q,O'}}$ is also trivial and 
$K_{X_{\q,O'}} = (\mu_{\q})^*K_{Y_{\l(\q),O'}}$. The second assertion follows from 
Proposition (3.3.3).  
\vspace{0.2cm} 

{\bf Corollary (3.3.6)}. {\em Assume that $$\nu^n: G \times^Q (\n(\q) + \tilde{O'}) \to 
\tilde{O}$$ is birational. Then, for $\q' \in \mathcal{S}(\l(\q))$, the normalized map 
$$(\nu')^n: G \times^{Q'}(\n(\q') + \tilde{O'}) \to \tilde{O}$$ is birational.} 
\vspace{0.15cm}

{\em Proof}. The map $$\mu_{\q',0}: X_{\q',O',0} \to Y_{\l(\q),O',0}$$ coincides 
with $(\nu')^n: G \times^{Q'}(\n(\q') + \tilde{O'}) \to \tilde{O}$ by (3.3.4). 
Moreover, $(\nu')^n$ is a generically finite morphism. 
On the other hand, $\mu_{\q'}$ is birational by (3.3.5). Since $Y_{\l(\q),O'}$ is 
normal, $\mu_{\q',0}$ has connected fibers. This menas that $(\nu')^n$ is birational. 
\vspace{0.2cm}

(3.4) {\em Nef cones and flops}: Let $(\q,O')$ be the same as in (3.3). 
Furthermore we assume that $$\nu^n: G \times^Q (\n(\q) + \tilde{O'}) \to \tilde{O}$$ 
is birational. 

(3.4.1) ({\em Nef cone of $G/Q$}): Let $Q \subset G$ be a parabolic subgroup 
of $G$.  Note that $\mathrm{Pic}(G/Q)  
\cong H^2(G/Q, \mathbf{Z})$.  
Define $M(L(Q)) := \mathrm{Hom}_{alg.gp}(L(Q), \mathbf{C}^*)$. 
Let $\chi: L(Q) \to \mathbf{C}^*$ be an element of $M(L(Q))$. 
By the exact sequence  
$$ 1 \to U(Q) \to Q \to L(Q) \to 1,$$ $\chi$ defines a group 
homomorphism $Q \to \mathbf{C}^*$, which gives rise to a 
line bundle $L_{\chi} := G \times^Q \mathbf{C}$ on $G/Q$.  
The correspondence $\chi \to L_{\chi}$ gives a map 
$$ \phi: M(L(Q)) \to \mathrm{Pic}(G/Q) $$ and it turns out 
be an isomorphism after tensorized with $\mathbf{R}$: 
$M(L(Q))_{\mathbf R} \cong \mathrm{Pic}(G/Q)_{\mathbf R}$.   
The nef cone $\overline{\mathrm{Amp}}(G/Q)$ is a closed convex cone in 
$H^2(G/Q, \mathbf{R})$ generated by nef line bundles on $G/Q$.  
Let us describe $\overline{\mathrm{Amp}}(G/Q)$ as a cone in 
$M(L(Q))_{\mathbf R}$ in terms of roots. Assume that $Q$ is a standard 
parabolic subgroups $Q_I$ containing a Borel subgroup $B$ (cf. (3.1)). 
Recall that $\Delta\setminus I$ corresponds to the set of marked vertices 
$\{v_1, ..., v_{\rho}\}$ of the Dynkin diagram.  Then  
$\rho = b_2(G/Q_I)$. 
The nef cone $\overline{\mathrm{Amp}}(G/Q_I)$ is then  
generated by dominant characters $\chi$ of $L(Q_I)$ (i.e. $\langle \chi, 
\alpha^{\vee}\rangle   
\geq 0, \forall \alpha \in \Delta$, where $\alpha^{\vee} \in \h$ is the coroot 
corresponding to $\alpha$, and $\chi$ is regarded as an element of 
$\h^*$). Moreover, it is a simplicial cone and its codimension 1 face 
consists of the dominant characters with $\langle \chi, v_i^{\vee}\rangle 
= 0$ for some $i$. We denote by $F_{v_i}$ this face.  
\vspace{0.2cm}
 
(3.4.2) ({\em Nef cones of $X_{\q, O'}$ and $X_{\q,O',0}$}):   
Note that each fiber of the natural projection $\pi: X_{\q,O'} \to G/Q$ is isomorphic 
to the normalization of $\r(\q) + \bar{O'}$, which coincides with $\r(\q) \times 
\tilde{O'}$. Here $\tilde{O'}$ is the normalization of $\tilde{O'}$. 
Since $\r(\q) \times \tilde{O'}$ is topologically contractible to the origin (by the 
natural $\mathbf{C}^*$-action), we see that $\pi^*: H^2(G/Q, \mathbf{R}) \cong 
H^2(X_{\q, O'}, \mathbf{R})$. Define $\overline{\mathrm{Amp}}(\mu_{\q})$ to be the 
closed convex cone in $H^2(X_{\q, O'}, \mathbf{R})$ generated by $\mu_{\q}$-nef 
line bundles on $X_{\q,O'}$. As in [Na 2], (P.3), one can prove that 
$$\pi^*(\overline{\mathrm{Amp}}(G/Q)) = \overline{\mathrm{Amp}}(\mu_{\q}).$$ 
Next let us consider $X_{\q,O',0}$. Let $p: X_{\q,O',0} \to G/Q$ be 
the natural projection. The situation is quite similar to the case of $X_{\q,O'}$. 
One can prove that 
$$ p^*(\overline{\mathrm{Amp}}(G/Q)) = \overline{\mathrm{Amp}}(\mu_{\q,0}).$$ 
In particular, $\overline{\mathrm{Amp}}(\mu_{\q})$ and 
$\overline{\mathrm{Amp}}(\mu_{\q,0})$ are both rational simplicial cones. 
Finally we shall define the {\em movable cones} $\overline{\mathrm{Mov}}(\mu_{\q})$ for $\mu_{\q}$ 
to be the closed convex cone of $H^2(X_{\q,O'}, \mathbf{R})$ generated 
by the $\mu_{\q}$-movable line bundles. Here a line bundle $L \in \mathrm{Pic}(X_{\q,O'})$ 
is called $\mu_{\q}$-movable if the support of 
$$\mathrm{Coker}[(\mu_{\q})^*(\mu_{\q})_*L \to L]$$ has codimension 
$\geq 2$. Denote by $\mathrm{Mov}(\mu_{\q})$ the interior of 
$\overline{\mathrm{Mov}}(\mu_{\q})$. 
Similarly we define $\overline{\mathrm{Mov}}(\mu_{\q,0})$ 
and $\mathrm{Mov}(\mu_{\q,0})$.    
\vspace{0.2cm}
   
(3.4.3) {\em Twists and flops}:     
We may assume that $Q$ is a standard parabolic subgroup $Q_I$ defined in (3.1).   
Take $v \in \Delta - I$ and put $\bar{I} := I \cup \{v\}$.  We shall use the same 
notation as in (3.2). Then  
$\bar{I}$ is decomposed into two sets $$\bar{I} = I_v \cup I'_v.$$ 
Let $\l({\bar I})$ be the standard Levi factor of $\q_{\bar I}$ and 
let $\z(\l({\bar I}))$ be its center. Then $\l({\bar I})/\z(\l({\bar I}))$ 
is decomposed into a simple factor $\l_{I_v}$ and other parts $\l_{I'_v}$: 
$$\l({\bar I})/\z(\l({\bar I})) = \l_{I_v} \oplus \l_{I'_v}.$$  
Note that $O' \subset \l_{I_v}$ 
or $O' \subset \l_{I'_v}$. For simplicity, we write $\g_v$ for $\l_{I_v}$ and 
$\bar{Q}$ for $Q_{\bar I}$.  Remark that $\q_v := \q \cap \g_v$ is a parabolic 
subalgebra of $\g_v$.   
Identify $\overline{\mathrm{Amp}}(\mu_{\q})$ with 
$\overline{\mathrm{Amp}}(G/Q)$ as in (3.4.2). As in (3.4.1), $v$ determines a codimension 1 face  
$F_v$ of $\overline{\mathrm{Amp}}(\mu_{\q_I})$.   
We are now going to construct the birational contraction map 
of $X_{\q,O'}$ corresponding to $F_v$.  
First look at the $\bar{Q}$-orbit $\bar{Q}\cdot(\r(\q) + \bar{O'})$ of 
$\r(\q) + \bar{O'}$. Then we can write 
$$\bar{Q}\cdot (\r(\q) + \bar{O'}) = \left\{ 
\begin{array}{rl}    
\r(\bar{\q}) \times G_v\cdot (\r(\q_v) + \bar{O'}) & \quad (O' \subset \g_v)\\ 
\r(\bar{\q}) \times \bar{O'} \times G_v \cdot \r(\q_v) & \quad (O' \subset \l_{I'_v}) 
\end{array}\right.$$
Put $W_v := W(\g_v)$ and let $\h_v$ be a Cartan subalgebra of $\g_v$.  
By the adjoint quotient map $G_v\cdot (\r(\q_v) + \bar{O'}) \to \h_v/W_v$ (or 
$G_v \cdot \r(\q_v) \to \h_v/W_v$), we have a map 
$\bar{Q}\cdot (\r(\q) + \bar{O'}) \to \h_v/W_v$. 
As in (3.3), we define a map $\eta: \bar{Q} \times^Q (\r(\q) + \bar{O'}) \to 
\k(\q_v)$. Since there is a natural map $\bar{Q}\times^Q (\r(\q) + \bar{O'}) \to  
\bar{Q}\cdot (\r(\q) + \bar{O'})$,  we get a map 
$$ \alpha: \bar{Q}\times^Q (\r(\q) + \bar{O'}) \to \bar{Q}\cdot (\r(\q) + \bar{O'}) 
\times_{\h_v/W_v} \k(\q_v).$$
Since $\alpha$ is a $\bar{Q}$-equivariant map, we obtain  
a map:  $$f : X'_{\q, O'} \to G \times^{\bar{Q}}{\bar{Q}}\cdot (\r(\q) + \bar{O'}) 
\times_{\h_v/W_v} \k(\q_v).$$ 
Note that $f$ is a morphism over $\k(\q)$. By restricting $f$ to the fibers over $0 \in \k(\q)$, 
we get a map $$f_0: G \times^Q (\n(\q) + \bar{O'}) \to G \times^{\bar Q}\bar{Q} \cdot (\n(\q) + \bar{O'}).$$ 
Define a map $f_v$ by $$f_v : \left\{ 
\begin{array}{rl}    
G_v \times^{Q_v} (\r(\q_v) + \bar{O'})  \to G_v \cdot (\r(\q_v) + \bar{O'}) 
\times_{\h_v/W_v} \k(\q_v) & \quad (O' \subset \g_v)\\ 
G_v \times^{Q_v} \r(\q_v) \to G_v \cdot \r(\q_v) 
\times_{\h_v/W_v} \k(\q_v) & \quad (O' \subset \l_{I'_v}) 
\end{array} \right.$$
Moreover define $\nu_v$ to be the generalized Springer map:  
$$\nu_v : \left\{ 
\begin{array}{rl}    
G_v \times^{Q_v} (\n(\q_v) + \bar{O'}) \to G_v \cdot (\n(\q_v) + \bar{O'}) & \quad (O' \subset \g_v)\\ 
G_v \times^{Q_v} \n(\q_v) \to G_v \cdot \n(\q_v) 
& \quad (O' \subset \l_{I'_v})
\end{array}\right.$$

Note that $$\bar{Q}\times^Q (\r(\q) + \bar{O'}) = 
\left\{ \begin{array}{rl}
\r(\bar{\q}) \times G_v \times^{Q_v} (\r(\q_v) + \bar{O'})
& \quad (O' \subset \g_v)\\ 
\r(\bar{\q}) \times \bar{O'} \times G_v \times^{Q_v}\r(\q_v)
& \quad (O' \subset \l_{I'_v}) 
\end{array}\right.,$$
$$\bar{Q}\times^Q (\n(\q) + \bar{O'}) = 
\left\{ \begin{array}{rl}
\n(\bar{\q}) \times G_v \times^{Q_v} (\n(\q_v) + \bar{O'})
& \quad (O' \subset \g_v)\\ 
\n(\bar{\q}) \times \bar{O'} \times G_v \times^{Q_v}\n(\q_v)
& \quad (O' \subset \l_{I'_v}) 
\end{array}\right.,$$
and  
$$\bar{Q}\cdot (\n(\q) + \bar{O'}) 
 = 
\left\{ \begin{array}{rl}
\n(\bar{\q}) \times G_v\cdot (\n(\q_v) + \bar{O'}) 
& \quad (O' \subset \g_v)\\ 
\n(\bar{\q}) \times \bar{O'} \times G_v\cdot \n(\q_v) 
& \quad (O' \subset \l_{I'_v}) 
\end{array}\right..$$ Therefore, we have the following lemma. 
\vspace{0.2cm}

{\bf Lemma (3.4.4)}. {\em (1) $$X'_{\q, O'} = \left\{ 
\begin{array}{rl}    
G \times^{\bar Q} (\r(\bar{\q}) \times 
G_v \times^{Q_v} (\r(\q_v) + \bar{O'})) & \quad (O' \subset \g_v)\\  
G \times^{\bar Q} (\r(\bar{\q}) \times \bar{O'} \times 
G_v \times^{Q_v} \r(\q_v)) 
& \quad (O' \subset \l_{I'_v}) 
\end{array}\right.$$
and 
$$ G \times^{\bar{Q}}{\bar{Q}}\cdot (\r(\q) + \bar{O'})
\times_{\h_v/W_v} \k(\q_v) = $$
$$\left\{ 
\begin{array}{rl}    
G \times^{\bar Q} (\r(\bar{\q}) \times 
G_v \cdot (\r(\q_v) + \bar{O'}) 
\times_{\h_v/W_v} \k(\q_v))  & \quad (O' \subset \g_v)\\  
G \times^{\bar Q} (\r(\bar{\q}) \times \bar{O'} \times 
G_v \cdot \r(\q_v)  
\times_{\h_v/W_v} \k(\q_v)) & \quad (O' \subset \l_{I'_v}). 
\end{array}\right.$$
When $O' \subset \g_v$, one has $f = id_G \times^{\bar Q} (id_{\r(\q)} \times f_v)$. 
When $O' \subset \l_{I'_v}$, one has $f = id_G \times^{\bar Q} (id_{\r(\q) \times \bar{O'}} \times f_v)$.}

{\em (2) $$G \times^Q (\n(\q) + \bar{O'}) = 
\left\{ 
\begin{array}{rl}    
G \times^{\bar Q} (\n(\bar{\q}) \times 
G_v \times^{Q_v} (\n(\q_v) + \bar{O'})) & \quad (O' \subset \g_v)\\
G \times^{\bar Q} (\n(\bar{\q}) \times \bar{O'} \times 
G_v \times^{Q_v} \n(\q_v))   
& \quad (O' \subset \l_{I'_v}) 
\end{array}\right.$$
and 
$$G \times^{\bar Q}\bar{Q} \cdot (\n(\q) + \bar{O'}) = 
\left\{ 
\begin{array}{rl}    
G \times^{\bar Q} (\n(\bar{\q}) 
\times G_v \cdot (\n(\q_v) + \bar{O'})) & \quad (O' \subset \g_v)\\  
G \times^{\bar Q} (\n(\bar{\q}) \times \bar{O'}  
\times G_v \cdot \n(\q_v))
& \quad (O' \subset \l_{I'_v}). 
\end{array}\right.$$
When $O' \subset \g_v$, one has $f_0 = id_G \times^{\bar Q} (id_{\n(\q)} \times \nu_v)$. 
When $O' \subset \l_{I'_v}$, one has $f_0 = id_G \times^{\bar Q} (id_{\n(\q) \times \bar{O'}} \times \nu_v)$.}   
\vspace{0.2cm}

{\bf Lemma (3.4.5)}. (1) {\em $f_0$ is birational.} 

(2) {\em $\nu_v$ is birational.}

(3) {\em $f_v$ is birational.} 

(4) {\em $f$ is birational.} 

{\em Proof}. (1): By the assumption, the generalized Springer map 
$G \times^Q (\n(\q) + \bar{O'}) \to \bar{O}$ is birational. 
Since $f_0$ factorizes this map, $f_0$ is birational. 

(2): Since $f_0$ is birational by (1), we see that $\nu_v$ is birational by (3.4.4), (2). 

(3): If $\nu_v$ is birational, then $f_v$ is birational by (3.3.3). 

(4): Since $f_v$ is birational by (3), $f$ is also birational by (3.4.4), (1). Q.E.D.  
\vspace{0.2cm}

Let $Z_v$ be the normalization of $G \times^{\bar{Q}}{\bar{Q}}\cdot (\r(\q) + \bar{O'}) 
\times_{\h_v/W_v} \k(\q_v)$.  
By Lemma (3.4.5), (4),  the map $f$ induces a birational morphism 
$$ f^n: X_{\q,O'} \to Z_v.$$ This map $f^n$ is the desired birational contraction 
map corresponding to $F_v$. 
   
We next let $\q'$ be the parabolic subalgebra obtained from $\q$ by the twist of 
$v$.  Then we have  
$$\k(\q_v) = \k(\q'_v).$$ 
Thus, there is a diagram of birational morphisms 
$$ X_{\q, O'} \stackrel{f^n}\to 
Z_v  
\stackrel{(f')^n}\longleftarrow X_{\q', O'},$$ and we have  
$$\overline{\mathrm{Amp}}(\mu_{\q}) \cap 
\overline{\mathrm{Amp}}(\mu_{\q'}) = F_v.$$

(3.5).  Let $O \subset \g$ be a nilpotent orbit. 
Assume that a parabolic subalgebra $\q_0$ of $\g$ and a nilpotent 
orbit $O' \subset \l(\q_0)$ give a {\bf Q}-factorialization  
$\nu^n: G \times^{Q_0} (\n(\q_0) + \tilde{O'}) \to \tilde{O}$. 
We put $\l := \l(\q_0)$ and denote by $L$ the corresponding 
Levi subgroup of $G$. Let us consider $Y_{\l, O'}$ defined in (3.3).  Note that, for 
$\q \in \mathcal{S}(\l)$, the nef cone $\overline{\mathrm{Amp}}(\mu_{\q})$ 
is regarded as a cone in $M(L)_{\mathbf R}$.  
\vspace{0.2cm}

{\bf Theorem (3.5.1)}. {\em For $\q \in \mathcal{S}(\l)$, the birational map 
$\mu_{\q}: X_{\q,O'} \to Y_{\l,O'}$ is a {\bf Q}-factorial terminalization and 
is an isomorphism in codimension one.  
Any {\bf Q}-factorial terminalization of $Y_{\l, O'}$ is obtained in this way. 
If $\q \ne \q'$, then $\mu_{\q}$ and $\mu_{\q'}$ give different 
{\bf Q}-factorial terminalizations. Moreover,} 
$$M(L)_{\mathbf R} = \cup_{\q \in \mathcal{S}(\l)} \overline{\mathrm{Amp}}(\mu_{\q}).$$   

{\em Proof}. We shall first prove that $\mu_{\q'}^{-1}\circ \mu_{\q}: X_{\q,O'} --\to 
X_{\q',O'}$ is not an isomorphism for $\q \ne \q'$. 
By Lemma (3.3.1), for $t \in (\k_0)^{reg}$ there is an isomorphism 
$$\rho_t: G \times^L (t + \bar{O'}) \cong G \times^Q (t + \n(\q) + \bar{O'})$$ 
defined by $\rho_t([g, t + y']) = [g, t+y']$. 
In a similar way, we have an isomorphism 
$$\rho'_t: G \times^L (t + \bar{O'}) \cong G \times^{Q'} (t + \n(\q') + \bar{O'}).$$ 
Note that $\rho'_t \circ (\rho_t)^{-1}$ coincides with 
$\mu_{\q',t}^{-1}\circ \mu_{\q,t}$. Assume that 
$\mu_{\q'}^{-1}\circ \mu_{\q}$ is an isomorphism. For $g \in G$ and $q \in Q - Q'$, 
let us consider two curves in $X_{\q,O'}$: 
$C_t := \{\rho_t([g,0])\}$ and $D_t := \{\rho_t([gq,0])\}$. 
If we let $t \to 0$, then we have 
$$ \lim_{t \to 0} C_t = \lim_{t \to 0} D_t.$$ 
Define $C'_t := \mu_{\q'}^{-1}\circ \mu_{\q}(C_t)$ and 
$D'_t := \mu_{\q'}^{-1}\circ \mu_{\q}(D_t)$. 
Note that, for $t \in (\k_0)^{reg}$, we have $C'_t = \{\rho'_t([g,0])\}$ 
and $D'_t = \{\rho'_t([g,q])\}$. 
Then $$ \lim_{t \to 0}C'_t = [g,0] \in G \times^{Q'}(\n(\q') + \bar{O'}),$$ 
and $$\lim_{t \to 0}D'_t = [gq,0] \in G \times^{Q'}(\n(\q') + \bar{O'}).$$ 
These two points should coincide. But, since $Q \ne Q'$, this is a contradiction. 
Therefore, $\mu_{\q'}^{-1}\circ \mu_{\q}$ is not an isomorphism.  
We next prove that any {\bf Q}-factorial terminalization $\mu: X \to Y_{\l,O'}$ 
is of the form $\mu_{\q}$. Fix a $\mu$-ample line bundle $L$ on $X$. 
Let $L^{(0)} \in \mathrm{Pic}(X_{\q_0,O'})$ be its proper transform. If $L^{(0)}$ is 
$\mu_{\q_0}$-nef, then $X = X_{\q_0,O'}$. Assume that $L^{(0)}$ is not $\mu_{\q_0}$-nef. 
There is an extremal ray $\mathbf{R}_+[z] \subset \overline{NE}(\mu_{\q_0})$ 
such that $(L^{(0)}, z) < 0$. Let $F \subset \overline{\mathrm{Amp}}(\mu_{\q_0})$ be 
the corresponding codimension one face. By (3.4) one can find $\q_1 \in 
\mathcal{S}(\l)$ such that $$\overline{\mathrm{Amp}}(\mu_{\q_0}) \cap \overline{\mathrm{Amp}}(\mu_{\q_1}) 
= F.$$ As constructed in (3.4), we then have a flop $X_{\q_0,O'} --\to X_{\q_1,O'}$. 
We let $L^{(1)} \in X_{\q_1, O'}$ be the proper transform of $L^{(0)}$ and 
repeat the same procedure. Thus, we get a sequence of flops 
$$ X_{\q_0,O'} --\to X_{\q_1,O'} --\to X_{\q_2,O'} --\to ... $$ 
But, since $\mathcal{S}(\l)$ is a finite set, this sequence must terminate by the 
discrepancy argument (cf. [Na 1, Theorem 6.1], [KMM, Proposition 5-1-11]). 
As a consequence, $X = X_{\q_k, O'}$ for some $k$.   
\vspace{0.2cm}

{\bf Proposition (3.5.2)}. {\em Let $Q \subset G$ be a maximal parabolic subgroup (i.e. 
$b_2(G/Q) = 1$) and 
let $O' \subset \l(\q)$ be a nilpotent orbit. Assume that  
$\nu^n: G \times^Q (\n(\q) + \tilde{O'}) \to \tilde{O}$ is a {\bf Q}-factorial terminalization 
which is an isomorphism in codimension one. 
Then $(\q, O')$ is a primitive pair (cf. (3.2.1)).} 

{\em Proof}. As in (3.4.3), we may assume that $Q = Q_I$. Since $Q$ is maximal, 
$\Delta - I = \{v\}$. Let $Q'$ be the parabolic subgroup twisted by $v$. 
As in (3.4), we have a birational map (over $\k(\q)$): 
$$ \gamma: X_{\q, O'} --\to X_{\q', O'}.$$ 
By restricting this diagram to the fibers over $0 \in \k(\q)$, we have a birational map 
$$ \gamma_0: G \times^Q (\n(\q) + \tilde{O'}) --\to G \times^{Q'}(\n(\q') + \tilde{O'}).$$ 
Let $p \in G \times^Q(\n(\q) + \tilde{O'})$ be a point such that $\nu^n$ is an isomorphism at $p$. 
Then $\gamma$ is an isomorphism at $p \in X_{\q,O'}$. Let $L \in \mathrm{Pic}(X_{\q,O'})$ 
be a $\mu_{\q}$-ample line bundle and denote by $\gamma_*(L) \in \mathrm{Pic}(X_{\q',O'})$ 
the proper transform of $L$ by $\gamma$. Then we have 
$$ \gamma_*(L)\vert_{G \times^{Q'} (\n(\q') + \tilde{O'})} = (\gamma_0)_*(L\vert_{G \times^Q (\n(q) + \tilde{O'})}).$$ 
Since $\mu_{\q}$ and $\mu_{\q'}$ are different by (3.5.1), the left hand side is not 
$(\nu')^n$-ample; hence the right hand side is not so. This means that $\gamma_0$ is not 
an isomorphism. Suppose $v$ is of the second kind; then $Q$ and $Q'$ are conjugate by an element 
$w \in W$.  
Since $\r(\q)$ and $\r(\q')$ are conjugate 
by $w$, $\r(\q) \cap \h$ and $\r(\q') \cap \h$ are also conjugate 
by $w$. Note that $\k(\q) := \r(\q) \cap \h = \r(\q') \cap \h$ and 
then $\l(\q) = \g^{\k(\q)}$ (cf. (3.1)). This means that $w$ sends $\l(\q)$ to $\l(\q)$; 
hence $w \in W'$. Since $W'$ stabilizes $O'$ by Lemma (3.3.2), 
two {\bf Q}-factorial terminalizations $\nu^n: G \times^Q(\n(\q) + \tilde{O'}) \to \tilde{O}$ 
and $(\nu')^n: G \times^{Q'} (\n(\q') + \tilde{O'}) \to \tilde{O}$ are the same one. In other words, 
$\gamma_0$ is an isomorphism; hence $v$ should be of the first kind. 
Let $D$ be the single marked Dynkin diagram corresponding to $Q$. Then $\l(\q)$ has only simple factors of type $A$ 
except when $D$ is of type $E_{6,I}$. If all simple factors of $\l(\q)$ are of type $A$, then 
$\tilde{O'}$ has {\bf Q}-factorial terminal singularities only when $O' = 0$. 
On the other hand, if $D$ is of type $E_{6,I}$, then $\l(\q) = D_5$. In $D_5$, we only have three 
nilpotent orbits $O'$ for which $\tilde{O'}$ has {\bf Q}-factorial terminal singularities: 
$O' = 0$, $O' = O_{[3,2^2,1^3]}$ and $O' = O_{[2^2,1^6]}$. 
\vspace{0.2cm}

{\bf Corollary (3.5.3)}. {\em Assume that $\nu^n: G \times^Q (\n(\q) + \tilde{O'}) \to \tilde{O}$  
is a {\bf Q}-factorial terminalization. 
Let $(f_0)^n: G \times^Q (\n(\q) + \tilde{O'}) \to Z_{v,0}$ be 
the birational contraction map corresponding to a codimension one face $F_v$ of $\overline{\mathrm{Amp}}(\mu_{\q})$ 
(cf.(3.4)). Then $(f_0)^n$ is an isomorphism in codimension one if and only if $(\q_v, O')$ is a primitive 
pair.}  

{\em Proof}. This follows from (3.4.4),(2) and (3.5.2). Q.E.D.  
\vspace{0.2cm}

Let us return to the original situation of (3.5). Define $\mathcal{S}^1(\l)$ to be the subset of 
$\mathcal{S}(\l)$ consisting of the parabolic subalgebras $\q$ obtained from $\q_0$ by a finite succession 
of the twists of the first kind. 
\vspace{0.2cm}

{\bf Theorem (3.5.4)}. {\em There is a one-to-one correspondence between the set of 
{\bf Q}-factorial terminalizations of $\tilde{O}$ and $\mathcal{S}^1(\l)$. 
In other words, every {\bf Q}-factorial terminalization of $\tilde{O}$ is obtained as  
$\mu_{\q,0}: X_{\q,O',0} \to \bar{O}$ for $\q \in 
\mathcal{S}^1(\l)$. Two different {\bf Q}-factorial terminalizations of 
$\tilde{O}$ are connected by a sequence of Mukai flops (cf. (3.2.1)). 
Moreover $$\overline{\mathrm{Mov}}(\mu_{\q_0,0}) = \cup_{\q \in 
\mathcal{S}^1(\l)} \overline{\mathrm{Amp}}(\mu_{\q,0}).$$} 

{\em Proof}. Let $\nu: Z \to \tilde{O}$ be a {\bf Q}-factorial terminalization of 
$\tilde{O}$. Note that the birational map $X_{\q_0,0} --\to Z$ is an isomorphism 
in codimension one. 
Fix a $\nu$-ample line bundle $M \in \mathrm{Pic}(Z)$. 
Let $M^{(0)} \in \mathrm{Pic}(X_{\q_0,O',0})$ be its proper transform. 
Assume that $M^{(0)}$ is not $\mu_{\q_0,0}$-nef. 
There is an extremal ray $\mathbf{R}_+[z] \subset \overline{NE}(\mu_{\q_0,0})$ 
such that $(M^{(0)}, z)  < 0$. Let $F \subset \overline{\mathrm{Amp}}(\mu_{\q_0,0})$ be 
the corresponding codimension one face. As explained in (3.4.3), $F$ corresponds 
to a marked vertex $v$ of the marked Dynkin diagram $D$ determined by $\q_0$. 
By Corollary (3.5.3) we see that $v$ is of the first kind (see (3.2)). 
By (3.4) one can find $\q_1 \in 
\mathcal{S}^1(\l)$ such that $$\overline{\mathrm{Amp}}(\mu_{\q_0,0}) \cap \overline{\mathrm{Amp}}(\mu_{\q_1,0}) 
= F.$$ As constructed in (3.4), we then have a flop $X_{\q_0,O',0} --\to X_{\q_1,O',0}$. 
We let $M^{(1)} \in X_{\q_1, O',0}$ be the proper transform of $M^{(0)}$ and 
repeat the same procedure. Thus, we get a sequence of flops 
$$ X_{\q_0,O',0} --\to X_{\q_1,O',0} --\to X_{\q_2,O',0} --\to ... $$ 
But, since $\mathcal{S}^1(\l)$ is a finite set, this sequence must terminate by the 
discrepancy argument (cf. [Na 1, Theorem 6.1], [KMM, Proposition 5-1-11]). 
As a consequence, $Z = X_{\q_k, O',0}$ for some $k$. 
\vspace{0.2cm}

(3.6) {\em Movable cones and the $W'$-action}:  
Start with the situation in Theorem (3.5.4). Recall that 
$$ W' := N_W(L_0)/W(L_0).$$ 
For $w \in N_W(L_0)$ and $\chi \in M(L_0)$, we define 
$w\chi  \in M(L_0)$ by $w\chi(g) = \chi (w^{-1}gw)$ with $g \in L_0$. 
In this way $N_W(L_0)$ acts on $M(L_0)_{\mathbf R}$. 
Note that $W(L_0)$ coincides with the subgroup of $N_W(L_0)$ 
which consists of the elements acting trivially on $M(L)_{\mathbf R}$. 
Hence, $W'$ acts on $M(L_0)_{\mathbf R}$ effectively. 
\vspace{0.2cm}

{\bf Theorem (3.6.1)}. 

{\em (i) The set $\mathcal{S}(\l)$ contains exactly $N \cdot \sharp(W')$ elements, where 
$N$ is the number of the conjugacy classes of parabolic subalgebras 
contained in $\mathcal{S}(\l)$.}
\vspace{0.2cm}

{\em (ii) For any $\q \in \mathcal{S}(\l)$, there is 
an element $w \in N_W(L_0)$ such that $w(\q) \in \mathcal{S}^1(\l)$.} 
\vspace{0.2cm}

{\em (iii) For any non-zero element $w \in W'$, we have 
$$ w(\mathrm{Mov}(\mu_{\q_0,0})) \cap \overline{\mathrm{Mov}}(\mu_{\q_0,0}) 
= \emptyset.$$} 

{\em (iv) The set $\mathcal{S}^1(\l)$ contains exactly $N$ elements.}   
\vspace{0.2cm}

{\em Proof}.  (i): 
Take two conjugate elements $\q, \q' \in \mathcal{S}(\l)$.  Then, 
there is an element $w \in W$ 
such that $\q = w(\q')$. Since $\r(\q)$ and $\r(\q')$ are conjugate 
by $w$, $\r(\q) \cap \h$ and $\r(\q') \cap \h$ are also conjugate 
by $w$. If we put $\k := \r(\q) \cap \h = \r(\q') \cap \h$, 
then $\l = \g^{\k}$ (cf. (3.1)). This means that $w$ sends $\l$ to $\l$; 
hence $w \in N_W(L)$.        
We next show that, if  $w(\q) = \q$ for an element $\q \in \mathcal{S}(\l)$,  then 
$w \in W(L)$.   
Let $U$ be the unipotent radical of $Q$. 
Then one can write $Q = U\cdot L$. Now we suppose that 
$w$ is represented by an 
element of the normalizer group $N_G(T)$.  Since $w(Q) = Q$ 
and $N_G(Q) = Q$, $w \in Q$.    
Let us write $w = u\cdot l$ 
with $u \in U$ and $l \in L$. By assumption,  
$w(L) = L$. This means that $u(L) = L$.  
Since any two Levi subgroups of $Q$ are conjugate by a 
unique element of $U$ (cf. [Bo], 14.19),  we have $u = 1$, 
which implies that $w \in W(L)$.   
 
(ii):  Let us assume that $\q_0$ is a standard parabolic 
subalgebra determined by a marked Dynkin diagram $D$. 
By the definition of twists, any  
$\q \in \mathcal{S}(\l)$ is conjugate to a standard 
parabolic subalgebra determined by a marked Dynkin diagram 
$D'$ which is equivalent to $D$ (cf.(3.2)). In order to get $D'$ from 
$D$, we only need the twists of the first kind. This means that 
there is an element $\q' \in \mathcal{S}^1(\l)$ such that $\q'$ is 
conjugate to the standard parabolic subalgebra determined by $D'$.    

(iii), (iv):  We shall prove that any two distinct elements 
of $\mathcal{S}^1(\l)$ are not conjugate to each other. 
Suppose that $\q, \q' \in 
\mathcal{S}^1(\l)$ are conjugate to each other. 
Let us consider the diagram 
$$ X_{{\q},O'} \stackrel{\mu_{\q}}\to Y_{\l,O'} 
\stackrel{\mu_{\q'}}\leftarrow X_{\q',O'}. $$ 
Restrict the diagram over $0 \in \k$ to get  
$$ X_{\q,O',0} \stackrel{\mu_{\q,0}}\to Y_{\l,0} 
\stackrel{\mu_{\q',0}}\leftarrow X_{\q',O', 0}. $$
By (3.3.4), this diagram coincides with 
$$ G \times^Q (\n(\q) + \tilde{O'}) \stackrel{\nu^n}\to 
\tilde{O} \stackrel{(\nu')^n}\leftarrow 
G \times^{Q'} (\n(\q') + \tilde{O'}).$$ 
Since $Q$ and $Q'$ are conjugate, we see that 
$\nu^n$ and $(\nu')^n$ give the same {\bf Q}-factorial 
terminalization of $\tilde{O}$ by (3.3.2). 
In other words, the birational map $\mu_{\q',0}^{-1} \circ \mu_{\q,0}$ 
is an isomorphism.  
We shall prove that $\mu_{\q'}^{-1}\circ \mu_{\q}$ is 
an isomorphism. Let $L$ be a $\mu_{\q}$-ample line 
bundle on $X_{\q,O'}$ and let  $L' \in \mathrm{Pic}(X_{\q',O'})$  
be the proper transform of $L$ by      
$\mu_{\q'}^{-1}\circ \mu_{\q}$.  
By Theorem (3.5.1),  
$X_{\q,O'}$ and $X_{\q', O'}$ are connected by a sequence of 
birational transformations which are isomorphisms in codimension 
one.  Since $\q, \q' \in \mathcal{S}^1(\l)$, these 
birational transformations all come from twists of the first kind. 
This means that, there is a closed subset $F$ of $X_{\q,O',0}$ with 
codimension $\geq 2$ such that $\mu_{\q'}^{-1}\circ \mu_{\q}$ 
is an isomorphism at each $x \in X_{\q,O',0}\setminus F$.  
Hence we have   
$$L'\vert_{X_{\q',O', 0}} \cong 
(\mu_{\q',0}^{-1}\circ \mu_{\q,0})_*(L\vert_{X_{\q,O',0}}).$$ 
But the right hand side is a $\mu_{\q',0}$-ample line bundle. 
Hence $L'\vert_{X_{\q',O',0}}$ is $\mu_{\q',0}$-ample. 
This shows that $L'$ is $\mu_{\q'}$-ample. Indeed, by the 
$\mathbf{C}^*$-action of $X_{\q',O'}$, every proper curve $C$ 
in a fiber of $\mu_{\q'}$ is deformed to a curve inside   
 $X_{\q',O',0}$; hence $(L', C) > 0$ follows from the 
ampleness of $L'\vert_{X_{\q',O',0}}$. Therefore,    
$\mu_{\q'}^{-1}\circ \mu_{\q}$ is 
an isomorphism.  Then, by Theorem (3.5.1),  $\q = \q'$. 
\vspace{0.2cm}

{\bf \S 4. Poisson deformations of nilpotent orbits. } 

(4.1) Let $X$ be a normal variety with symplectic singularities (cf. (2.4)). 
We shall define a Poisson structure on $X$ by using the symplectic 2-form $\omega$ on $X_{reg}$. 
By $\omega$ the sheaf of 1-forms $\Omega^1_{X_{reg}}$ is identified with the 
sheaf of vector fields $\Theta_{X_{reg}}$. Since $\Omega^2_{X_{reg}} \cong \wedge^2 \Theta_{X_{reg}}$, 
$\omega$ determines a bivector $\Theta \in \wedge^2 \Theta_{X_{reg}}$. 
We then define a bracket 
$$ \{\:, \:\} : \wedge^2_{\mathbf C} \mathcal{O}_{X_{reg}} \to \mathcal{O}_{X_{reg}} $$ 
by $\{f,g\} := \Theta (df \wedge dg)$. By definition this bracket is bi-derivation. 
Moreover, it satisfies the {\em Jacobi identity} 
$$ \{\{f,g\},h\} + \{\{g, h\}, f\} + \{\{h,f\}, h\} = 0 $$  
(cf. [C-G], Theorem 1-2-7). In other words, we have a Poisson structure on 
$X_{reg}$. Since $X$ is normal, this bracket uniquely extends to the bracket 
$$ \{\:, \:\} : \wedge^2_{\mathbf C} \mathcal{O}_{X} \to \mathcal{O}_{X}. $$
This bracket is also a bi-derivation and satisfies the Jacobi-identity. 
In this way, $(X, \{\:, \:\})$ is a variety with a Poisson structure.    
We shall introduce the notion of a Poisson deformation of $(X, \{\:, \:\})$.  
First recall    
\vspace{0.2cm} 

{\bf Definition (4.1.1)}. Let $T$ be a scheme (resp. complex space).  
Let $\mathcal{X}$ be a scheme  (resp.  
complex space) over $T$. 
Then $(\mathcal{X}, \{\; , \;\})$ is a Poisson scheme (resp. a Poisson space) over $T$ if 
$\{\; , \;\}$ is an $\mathcal{O}_T$-linear map: 
$$ \{\; , \;\}: \wedge^2_{\mathcal{O}_T}\mathcal{O}_{\mathcal{X}} \to 
\mathcal{O}_{\mathcal{X}}$$ 
such that, for $a,b,c \in \mathcal{O}_{\mathcal{X}}$, 
\begin{enumerate}
\item $\{a, \{b,c\}\} + \{b,\{c,a\}\} + \{c,\{a,b\}\} = 0$ 
\item $\{a,bc\} = \{a,b\}c + \{a,c\}b.$ 
\end{enumerate}

Let $0 \in T$ be a punctured $\mathbf{C}$-scheme.
A Poisson deformation of $(X, \{\;, \;\})$ 
over $T$ is a pair of a Poisson scheme $(\mathcal{X}, \{\;, \;\}_T)$ over $T$ 
and an isomorphism $\phi: \mathcal{X} \times_T \mathrm{Spec}(\mathbf{C}) 
\cong X$ 
such that $\mathcal{X}$ is flat over $T$, 
and the Poisson structure $\{\; , \;\}_T$ induces 
the original Poisson structure $\{\; , \;\}$ over the closed 
fiber $X$ by $\phi$.   
Let $S$ be a local Artin $\mathbf{C}$-algebra with 
residue field $\mathbf{C}$. Two Poisson deformations $(\mathcal{X}, \phi)$ and 
$(\mathcal{X}', \phi')$ over $S$ are {\em equivalent} if there is a Poisson 
isomorphism $\varphi: \mathcal{X} \cong \mathcal{X}'$ over 
$\mathrm{Spec}(S)$ which induces the identity map of $X$ over 
$\mathrm{Spec}(\mathbf{C})$ via $\phi$ and $\phi'$. 

The Poisson deformation $\mathcal{X} \stackrel{f}\to T$ 
is called {\em formally universal}  at $0 \in T$ if, for any 
Poisson deformation $\mathcal{X}' \to T'$ of $X$ with a 
local Artinian base $T'$, there is a unique map $T' \to T$ 
such that $\mathcal{X}' \cong \mathcal{X}\times_T T'$ as 
a Poisson deformation of $X$ over $T'$. 
In this case, for a small open neighborhood $V$ of $0 \in T$, the family $f\vert_{f^{-1}(V)}: 
f^{-1}(V) \to V$ is called the Kuranishi family for the Poisson deformations of $X$, and 
$V$ is called the Kuranishi space for the Poisson deformations of $X$. 

{\bf Proposition (4.2)}. {\em Let $O$ be a nilpotent orbit of a complex simple Lie algebra $\g$ 
and let $\nu^n: G \times^Q (\n(\q) + \tilde{O}') \to \tilde{O}$ be a 
{\bf Q}-factorial terminalization. Then $G \times^Q (\n + \tilde{O}')$ has symplectic 
singularities. Moreover,  
$$ X_{\q, O'} := G \times^Q (\r(\q) + \tilde{O}') \to \k(\q) $$ is a 
Poisson deformation of  $G \times^Q (\n(\q) + \tilde{O}')$.}  

{\em Proof}.  For $t \in \k(\q)$, the fiber $X_{\q, O',t}$ is isomorphic 
to $G \times^Q (t + \n(\q) + \tilde{O}')$, whose regular locus 
is $G \times^Q (t + \n(\q) + O')$. 
We have a natural $G$-equivariant map 
$$\mu_t: G \times^Q (t + \n(\q) + \tilde{O}') \to \g$$ 
defined by $[g, t + y + y'] \to Ad_g(t+y+y')$. 
The image of this map coincides with the closure of an 
adjoint orbit, say $O_{\mu_t}$. We shall prove that the 
pull-back $\omega_t$ of the Kostant-Kirillov 2-form on $O_{\mu_t}$ 
give a symplectic 2-form on 
$G \times^Q (t + \n(\q) + O')$.  
Let $\l$ be the Levi part of $\p$ and 
fix a Cartan subalgebra $\h$ of $\g$ such that 
$\h \subset \l$. In the remainder of the proof, we shall simply write $\n$ for 
$\n(\q)$. 
There is an involution $\phi_{\g}$ of $\g$ which stabilizes 
$\h$ and which acts on the root system $\Phi$ via $-1$. Put $\n_{-} := \phi_{\g}(\n)$. 
Take a point $[1, t+y+y'] \in G \times^Q (t + \n + {\mathcal{O}}')$ so that 
$y \in \n(\q)$, $y' \in \mathcal{O}'$ 
and $t + y + y' \in O_{\mu_t}$. The tangent space of 
$G \times^Q (t + \n + \bar{O}')$ at $[1, t+y+y']$ is decomposed as 
$$ T_{[1, t+y+y']} = \n_{-} \oplus T_{y+y'}(t + \n + \bar{O}'). $$ 
Since $Q\cdot (t+y+y')$ coincides with the Zariski open dense subset 
$O_{\mu_t} \cap (t + \n + {O}') \subset t + \n + {O}'$,  
an element $v \in T_{[1, t+y+y']}$ can be written as 
$$ v = v_1 + [v_2, t+y+y'], \;\; v_1 \in \n_{-}, \; v_2 \in \q.$$ 
Let $d\nu_*: T_{[1, y+y']} \to T_{\nu([1, t+y+y'])}O_{\mu_t}$ be the tangential map 
for $\mu_t$. 
Then $$d(\mu_t)_*(v) = [v_1 + v_2, t+y+y'].$$ 
Take one more element $w \in T_{[1, t+y+y']}$ in such a way that 
$$ w = w_1 + [w_2, t+y+y'], \;\; w_1 \in \n_{-}, \; w_2 \in \q.$$ 
Denote by $\langle \;, \;\rangle$ the Killing form of $\g$. 
By the definition of the Kostant-Kirillov form, one has 
$$ \omega(d(\mu_t)_*(v), d(\mu_t)_*(w)) := \langle t+y+y', [v_1 + v_2, w_1 + w_2]\rangle.$$ 
Note that $\langle t+y+y', [v_1, w_1]\rangle = \langle y, [v_1, w_1]\rangle$, and 
$\langle t+y+y', [v_2, w_2]\rangle = \langle y', [v_2, w_2]\rangle$.  Therefore, 
\vspace{0.2cm}

$\omega(d(\mu_t)_*(v), d(\mu_t)_*(w)) =  $ 

$\langle y, [v_1,w_1] \rangle +  \langle t+y+y', [v_1, w_2] \rangle   
+ \langle t+y+y', [v_2, w_1] \rangle + \langle y', [v_2, w_2] \rangle = $ 

$\langle y, [v_1,w_1] \rangle + \langle v_1, [w_2, t+y + y']_n\rangle  - 
\langle [v_2, t+y + y']_n, w_1 \rangle  
+ \omega ([v_2, y'], [w_2, y']), $
\vspace{0.2cm}

where 
$[w_2, t+y + y']_n$ (resp. $[v_2, t+y + y']_n$) is 
the nil-radical part of  $[w_2, t+y + y']$ (resp. 
$[v_2, t+y + y']$) in the decomposition 
$T_{y+y'}(t+ \n + \bar{\mathcal O}') = \n + T_{y'}\mathcal{O}'$. 
Let $\mathcal{O}_r \subset \g$ be the Richardson orbit for $Q$, and 
let $\pi: T^*(G/Q) \to \bar{\mathcal{O}_r}$ be the Springer map. 
The first part $\langle y, [v_1,w_1] \rangle + \langle v_1, [w_2, t+y + y']_n\rangle  - 
\langle [v_2, t+y + y']_n, w_1 \rangle$ 
corresponds to the 2-form on $T^*(G/Q)$ obtained by the pull-back of 
the Kostant-Kirillov 2-form on $\mathcal{O}_r$ by $\pi$ (cf. [Pa]), which is non-degenerate 
on $T^*(G/Q)$. 
Let us consider the second part $\omega ([v_2, y'], [w_2, y'])$. 
Denote by $[v_2, t+y+y']_l$ (resp. 
$[w_2, t+y+y']_l$) the $T_{y'}\mathcal{O}'$-part of $[v_2, t+y+y']$ (resp. $[w_2, t+y+y']$) in 
the decomposition 
$T_{t+y+y'}(t+ \n + \bar{\mathcal O}') = \n + T_{y'}\mathcal{O}'$. 
Then 
$[v_2,y'] = [v_2, t+y+y']_l$ and $[w_2, y'] = [w_2, t+y+y']_l$; hence,  
the second part is the Kostant-Kirillov form on $\mathcal{O}'$. 

Now assume that $y \in \n$ and $y' \in {\mathcal O}'$ (not necessarily 
$t + y + y' \in O_{\mu_t}$).  Write $v \in T_{[1, t+y+y']}$ as 
$$ v = v_1 + (v_3)_n + (v_3)_l, $$ where $v_1 \in \n_{-}$, $(v_3)_n \in \n$, 
and $(v_3)_l \in T_{y'}{\mathcal O}'$. Similarly, write $w \in T_{[1,y+y']}$ as 
$$ w = w_1 + (w_3)_n + (w_3)_l. $$ 
The arguments above show that      
$$(\mu_t)^*\omega(v,w) = \langle y, [v_1,w_1] \rangle + \langle v_1, (w_3)_n \rangle 
- \langle (v_3)_n, w_1 \rangle + \omega((v_3)_l, (w_3)_l).$$ It is easily checked that  
$(\mu_t)^*\omega$ is non-degenerate at $[1, t+y+y']$. By the $G$-equivariance 
of $\mu_t$, we see that $\omega_t := (\mu_t)^*\omega$ is a symplectic 2-form on 
$G \times^Q (t + \n(\q) + O')$.

By this description we can also observe that the family of symplectic 2-forms $\{\omega_t\}$ 
defines a relative symplectic 2-form $\omega$ of $G \times^Q (\r(\q) + O') \to \k(\q)$. 
This relative symplectic 2-form makes $G \times^Q (\r(\q) + \tilde{O}')$ into a 
Poisson scheme over $\k(\q)$. Its central fiber is clearly the original 
Poisson scheme $G \times^Q (\n + \tilde{O}')$. Q.E.D.  
\vspace{0.15cm}

(4.3) We have constructed in (3.3) a map $\chi: G \cdot (\r(\q) + \bar{O}') \to \h/W$ and  
have remarked that the normalization of $\mathrm{Im}(\chi)$ coincides with $\k(\q)/W'$.
Let $G \cdot (\r(\q) + \bar{O}')^n$ be the normalization of $G \cdot (\r(\q) + \bar{O}')$. 
Then $\chi$ induces a map $$\chi^n: G \cdot (\r(\q) + \bar{O}')^n \to \k(\q)/W'.$$ 

{\bf Proposition (4.3.1)}. {\em $\chi^n$ is a flat morphism whose central fiber is isomorphic to $\tilde{O}$.  
Moreover, $\chi^n: G \cdot (\r(\q) + \bar{O}')^n \to \k(\q)$ is a Poisson deformation 
of $\tilde{O}$.} 

{\em Proof}. When $O' = 0$, the statements are exactly Corollary (2.3) and Proposition (2.6) of 
[Na 5]. The proof in a general case is the same.  
\vspace{0.2cm}

(4.4) {\bf The period map}: We put $X := G \times^Q (\n(\q) + \tilde{O}')$ and 
consider the $\mathbf{C}^*$-equivariant Poisson deformation of $X$: 
$$ X_{\q, O'}  \to \k(\q). $$ 
In (4.2) we have defined a relative symplectic 2-form $\omega$ 
on the regular locus $(X_{\q, O'})_{reg} := G \times^Q(\r(\q) + O')$ over $\k(\p)$. 
We shall construct a period map 
$$ p: \k(\q) \to H^2(X, \mathbf{C}) $$ by using $\omega$. 
Since the fibers of $ X_{\q, O'}  \to \k(\q)$ are not smooth, we need 
some technical arguments to define $p$. First of all, note that $X$ is a {\bf Q}-factorial 
terminalization of $\tilde{O}$, where $\tilde{O}$ has a $\mathbf{C}^*$-action with 
positive weights. The $\mathbf{C}^*$-action on $X$ is the lifting of this 
$\mathbf{C}^*$-action. Then $X$ is also {\bf Q}-factorial as a complex analytic 
space by [Na 4], Proposition (A.9). By [ibid, Theorem 17] we see that 
$X_{\q, O'}  \to \k(\q)$ is a locally trivial flat deformation of $X$. 
In the proof of [ibid], Proposition 24, we have constructed a simultaneous 
$\mathbf{C}^*$-equivariant resolution of $X_{\q, O'}  \to \k(\q)$: 
$$\beta: \mathcal{Z} \to X_{\q,O'}.$$ 

We now have a commutative diagram 
\begin{equation} 
\begin{CD} 
\mathcal{Z} @>>> G \cdot(\r(\q) + \bar{O}')^n \\ 
@V{\alpha}VV @VVV \\ 
\k(\q) @>>> \k(\q)/W'.     
\end{CD} 
\end{equation} 

Note that $G \cdot(\r(\q) + \bar{O}')^n$ has a $\mathbf{C}^*$-action with a unique 
fixed point and with positive weights.
Moreover, the diagram is $\mathbf{C}^*$-equivariant and 
$\alpha: \mathcal{Z} \to \k(\q)$ 
is a simultaneous resolution of $G \cdot (\r(\q) + \bar{O}')^n \to \k(\q)/W'$. Then we see that 
$\mathcal{Z}$ is a $C^{\infty}$-trivial fiber bundle over $\k(\q)$  
by [Slo], Remark at the end of section 4.2. 
Let  $\Omega^{\cdot}_{\mathcal{Z}^{an}/\k(\q)}$ be the relative complex-analytic de Rham 
complex. Let $\mathcal{K}$ be the subsheaf of $\Omega^2_{\mathcal{Z}^{an}/\k(\q)}$ 
which consists of d-closed relative 2-forms. By the natural map $\mathcal{K}[-2] \to 
\Omega^{\cdot}_{\mathcal{Z}^{an}/\k(\q)}$, we can define a sequence of maps: 
$$ \alpha_*\mathcal{K} \to  \mathbf{R}^2\alpha_*\Omega^{\cdot}_{\mathcal{Z}^{an}/\k(\q)} 
\cong R^2\alpha_*\alpha^{-1}\mathcal{O}^{an}_{\k(\q)}.$$ 
Since $R^2\alpha_*\alpha^{-1}\mathcal{O}^{an}_{\k(\q)} \cong R^2\alpha_*\mathbf{C} \otimes_{\mathbf C}
\mathcal{O}^{an}_{\k(\q)}$ (cf. [Lo], Lemma (8.2)), we have an isomorphism 
$$R^2\alpha_*\alpha^{-1}\mathcal{O}^{an}_{\k(\q)} \cong H^2(\mathcal{Z}_0, \mathbf{C}) \otimes_{\mathbf C} 
\mathcal{O}^{an}_{\k(\q)}.$$ By pulling back the relative symplectic 2-form $\omega$ 
of $X_{\q,O'}/\k(\q)$ defined in (4,2), we get a relative d-closed 2-form $\omega_{\mathcal Z}$ 
of $\mathcal{Z}/\k(\q)$. 
Then $\omega_{\mathcal Z}$ gives a section $s$ of the 
sheaf $H^2(\mathcal{Z}_0, \mathbf{C}) \otimes_{\mathbf C} 
\mathcal{O}^{an}_{\k(\q)}$.  Let $$ev_t: H^2(\mathcal{Z}_0, \mathbf{C}) \otimes_{\mathbf C} 
\mathcal{O}^{an}_{\k(\q)} \to H^2(\mathcal{Z}_0, \mathbf{C})$$ be the evaluation map at 
$t \in \k(\q)$. We define a period map 
$$ p: \k(\q) \to H^2(\mathcal{Z}_0, \mathbf{C}) $$ by $p(t) = ev_t(s)$. By the construction, 
$p$ is a holomorphic map.
The birational morphism $\mathcal{Z} \to X_{\q,O'}$ induces a birational morphism 
$\mathcal{Z}_0 \to X_{\q,O',0}$ of central fibers. 
Since $X_{\q,O',0}$ has only rational singularities, there is an injection 
$H^2(X_{\q,O',0}, \mathbf{C}) \to H^2(\mathcal{Z}_0, \mathbf{C})$. 
We shall prove that $p$ factors through $H^2(X_{\q,O',0}, \mathbf{C})$: 
$$ p: \k(\q) \to H^2(X_{\q,O',0}, \mathbf{C}) \to H^2(\mathcal{Z}_0, \mathbf{C}).$$ 
Take a point $t \in \k(\q)$. In order to prove that 
$p(t) \in H^2(X_{\q,O',0}, \mathbf{C})$, it is enough to show that $(p(t), C) = 0$ 
for all proper curves $C$ which are mapped to points by the map $\beta_0: \mathcal{Z}_0 
\to X_{\q,O',0}$. Put $p := \beta_0(C) \in X_{\q,O',0}$ and take a small open 
neighborhood $V$ of $p \in X_{\q,O'}$. Set $\tilde{V} := \beta^{-1}(V)$. 
Let us consider the map $V \to \k(\q)$ obtained as the composite $V \subset X_{\q,O'} \to \k(\q)$. 
We regard this map as the flat deformation of the central fiber $V_0$. 
Let $L$ be the line of $\k(\q)$ passing through $0$ and $t$. 
Restrict the map $V \times_{\k(\q)} L \to L$ to the $n$-th infinitesimal neighborhoods $L_n$ of 
$0 \in L$. Then we have a formal deformation $\{V_n\}$ of $V_0$.   
The map $\beta$ induces a resolution $\tilde{V}_n \to V_n$ for each $n \geq 1$. 
Since $X_{\q,O'} \to \k(\q)$ is a locally trivial flat deformation of $X_{\q,O',0}$, 
we see that $V_n = V_0 \times L_n$. By the construction of $\beta$ (cf. [Na 4], 
Proposition 24), we have $\tilde{V}_n = \tilde{V}_0 \times L_n$. 
This means that the proper curve $C \subset \tilde{V}_0$ deforms sideways in 
the flat deformation $\tilde{V}_n \to L_n$ for each $n$. Here let us 
consider the relative Hilbert scheme 
$$H:= \mathrm{Hilb}(\mathcal{Z}\times_{\k(\q)}L/X_{\q,O'}\times_{\k(\q)}L).$$ 
The argument above shows that there is an irreducible component $H'$ of $H$ 
such that $[C] \in H'$ and $H'$ dominates $L$ by the composite 
$H' \to X_{\q,O'} \times_{\k(\q)}L \to L$. Note that $\mathcal{Z} \to \k(\q)$ 
is $\mathbf{C}^*$-equivariant and $L - \{0\}$ is a $\mathbf{C}^*$-orbit of 
$\k(\q)$. Since $H'$ dominates $L$, one can find $t' \in L - \{0\}$ such that 
there is a proper curve $C_{t'} \subset \mathcal{Z}_{t'}$ which is deformation 
equivalent to $C$. By using the $\mathbf{C}^*$-action, one can also find a 
proper curve $C_t \subset \mathcal{Z}_t$ which is deformation equivalent 
to $C$. By the definition of $p$, we have 
$$(p(t),C)_{\mathcal{Z}_0} = ([\omega_{\mathcal{Z}_t}], C_t)_{\mathcal{Z}_t}.$$ 
Since the restriction of a holomorphic 2-form to a curve is always zero, 
$([\omega_{\mathcal{Z}_t}], C_t)_{\mathcal{Z}_t} = 0$. 
\vspace{0.2cm}

Once we have made the period map precise, we can generalize Proposition (2.7) 
and Theorem (2.8) of [Na 5] as follows. The proof is almost the same as [Na 5]. 
  
{\bf Theorem (4.5)}. {\em Let $O$ be a nilpotent orbit of a 
complex simple Lie algebra $\g$, and let $\nu^n: G \times^Q (\n(\q) + \tilde{O}') \to \tilde{O}$
be a {\bf Q}-factorial terminaization. Then the following ${\mathbf C}^*$-equivariant commutative diagram} 

\begin{equation} 
\begin{CD} 
X_{\q,O'} @>>> G \cdot(\r(\q) + \bar{O}')^n \\ 
@VVV @VVV \\ 
\k(\q) @>>> \k(\q)/W'     
\end{CD} 
\end{equation} 
{\em gives formally universal deformations of  $G \times^Q (\n(\q) + \tilde{O}')$ 
and $\tilde{O}$.}

Yoshinori Namikawa 

Department of Mathematics, Faculty of Science, Kyoto University, 
Kyoto, 606-8502, Japan 

e-mail: namikawa@math.kyoto-u.ac.jp 

\end{document}